\documentclass[11pt]{amsart}
\usepackage[babel]{csquotes}
\usepackage{esint}
\usepackage{enumitem}
\usepackage{amsmath,amsthm,amssymb,mathrsfs,amsfonts,verbatim,enumitem,color,leftidx,mathtools}
\usepackage{mathabx}
\usepackage{etoolbox} 
\usepackage{bbm}
\usepackage[all,tips]{xy}
\usepackage{graphicx,ifpdf}
\usepackage{stmaryrd}
\usepackage{amssymb}
\usepackage{marginnote}
\usepackage{dsfont}

\ifpdf
\fi
\usepackage[colorlinks]{hyperref}
\hypersetup{
linkcolor=blue,          % color of internal links (change box color with linkbordercolor)
citecolor=green,        % color of links to bibliography
}

\usepackage{tikz}%pictures drawn in latex

\newtheorem{thm}{Theorem}[section]
\newtheorem{lem}[thm]{Lemma}

\newtheorem{prop}[thm]{Proposition}

\newtheorem{conj}[thm]{Conjecture}

\theoremstyle{definition}

\newtheorem{remark}[thm]{Remark}
\newtheorem{example}[thm]{Example}

\theoremstyle{remark}

\numberwithin{equation}{section}
\definecolor{esperance}{rgb}{0.0,0.5,0.0}

\newcommand{\bm}{\mathbf{m}}

\newcommand{\bp}{\mathbf{p}}
\newcommand{\bq}{\mathbf{q}}

\newcommand{\bx}{\mathbf{x}}

\newcommand{\bxi}{\boldsymbol{\xi}}

\newcommand{\R}{\mathbb{R}}
\newcommand{\N}{\mathbb{N}}
\newcommand{\Z}{\mathbb{Z}}
\newcommand{\Q}{\mathbb{Q}}
\newcommand{\C}{\mathbb{C}}

\newcommand{\bP}{\mathbb{P}}

\newcommand{\onto}{\xymatrix{\ar@{>>}[r]&}}
\makeatletter
\newcommand*{\rom}[1]{\expandafter\@slowromancap\romannumeral #1@}
\makeatother

\begin{document}

\title[On systems of quadratic forms in many variables]{On systems of quadratic forms in many variables}
\author{Wooyeon Kim and Kiseok Yeon}

% \thanks will become a 1st page footnote.
\thanks{}

%\author{}


\def\thefootnote{}
   %%d 
\def\thefootnote{\arabic{footnote}}
\setcounter{footnote}{0}

\begin{abstract}
In this paper, we prove that a system of $R$ quadratic forms in at least $4R+1$ variables has a nonzero integer solution, provided that it defines a smooth complete intersection of codimension $R$ over $\mathbb{C}$ and is soluble over $\mathbb{R}$ and $\mathbb{Q}_p$ for all primes $p\leq (2R)^R$. Furthermore, we establish an asymptotic formula for the number of integer solutions of bounded height to such a system. More precisely, let $Q_1,\ldots,Q_R\in \mathbb{Z}[\mathbf{x}]$ be quadratic forms in $n$ variables, and define
\begin{equation*}
    N(T):=\#\{\mathbf{x}\in [-T,T]^n\cap \mathbb{Z}^n:\   Q_1(\mathbf{x})=\cdots= Q_R(\mathbf{x})=0\}.
\end{equation*} 
Let $\mathcal Z_{\mathbb{C}}\subset\bP^{n-1}(\mathbb{C})$ be the
projective variety defined by $Q_1=\cdots=Q_R=0$, and suppose that $\mathcal Z_{\mathbb{C}}$ is a smooth complete intersection of codimension $R$. Then, whenever $n\geq 4R+1$, one has the asymptotic formula
\begin{equation*}
     N(T)=cT^{n-2R}+O(T^{n-2R-\delta}),\ \text{as}\ T\rightarrow\infty,
\end{equation*}
 for some $\delta>0,$ where $c$ is the product of local densities. The variable bound $n\geq  4R+1$ represents the square-root cancellation threshold.
\end{abstract}

\maketitle

\section{Introduction}
In this paper, we study integer solutions to a system of quadratic forms defining a smooth complete intersection. 
Let $R\geq1$, and let $Q_1,\ldots,Q_R\in\Z[x_1,\ldots,x_n]$ be quadratic forms with integer coefficients.
Consider the following system of quadratic equations:
\begin{equation}\label{systemofquadraticforms}
Q_i(x_1,\ldots,x_n)=0 \quad \textrm{ for all }1\leq i\leq R.
\end{equation}
Let $\mathcal Z\subset\bP^{n-1}_{\mathbb{Q}}$ be the
projective variety defined by $Q_1=\cdots=Q_R=0$.
For a field extension $K/\mathbb Q$, let $\mathcal Z_K$
denote its base change to $K$. Throughout this paper,
we assume that $\mathcal Z_\C$ is a smooth complete
intersection of codimension $R$. Explicitly, we require
\begin{equation}\label{eq:geometric-hypothesis}
\operatorname{rank}_{\C}
\bigl(\nabla Q_1(\bx),\ldots,\nabla Q_R(\bx)\bigr)=R
\end{equation}
for every $\bx=(x_1,\ldots,x_n)\in \mathcal{Z}_{\mathbb{C}}$.

We begin with a conjecture on integer solutions to the system (\ref{systemofquadraticforms}). It follows from the conjectural
Hasse principle for smooth complete intersections of $R$ quadrics
in more than $4R$ variables, discussed in \cite[p.367]{BHB17} (see also \cite[p.207]{RM18}),
together with the conjecture that every system of $R$ quadratic
forms in more than $4R$ variables has a nontrivial common zero
over $\mathbb{Q}_p$ for every prime $p$; see \cite{HB13}.
\begin{conj}\label{conj:Meyerforsystem}
Let $n\geq 4R+1$. Assume that $\mathcal{Z}_{\mathbb{C}}$ is a smooth complete intersection and $\mathcal{Z}_{\mathbb{R}}$ is nonempty. Then the system \eqref{systemofquadraticforms} has a nonzero solution $\mathbf{x}\in\mathbb{Z}^n$.
\end{conj}

The bound on the number of variables $n\geq4R+1$ in
Conjecture~\ref{conj:Meyerforsystem} is sharp.
Indeed, for every $R\geq1$ and every prime $p\equiv3\pmod4$,
there exist $R$ diagonal quadratic forms with integer
coefficients in $4R$ variables defining a smooth complete
intersection with a real point but no $\mathbb{Q}_p$-point,
and hence no rational point; see Example~\ref{example:sharpness}. The smoothness hypothesis is also essential: even for $R=2$
and arbitrarily large $n$, a complete intersection may have a
nonsingular real point but no rational point; see
Example~\ref{example:smoothness}.

\bigskip

Our first theorem provides a positive answer to Conjecture~\ref{conj:Meyerforsystem} for any $R\geq 1$, under the additional assumption of local solubility at small primes. For large primes $p$, Heath-Brown \cite{HB10} proved that the system $(\ref{systemofquadraticforms})$ is soluble over $\mathbb{Q}_p$ for $p>(2R)^R$.

\begin{thm}\label{thm:Meyerforsystem}
    Let $n\geq 4R+1$. Assume that $\mathcal{Z}_{\mathbb{C}}$ is a smooth complete intersection and $\mathcal{Z}_{\mathbb{R}}$ is nonempty. Suppose further that for all $p\leq (2R)^R$ the system \eqref{systemofquadraticforms} has a nonzero solution over $\mathbb{Z}_p$. Then the system \eqref{systemofquadraticforms} has a nonzero solution $\mathbf{x}\in\mathbb{Z}^n$.
\end{thm}

\bigskip

We now turn to quantitative results on integer solutions to the system \eqref{systemofquadraticforms}. The next theorem provides the asymptotic formula for the number of such solutions of bounded height. Define
\begin{equation*}
     N(T):=\#\{\mathbf{x}\in [-T,T]^n\cap \mathbb{Z}^n:\   Q_1(\mathbf{x})=\cdots= Q_R(\mathbf{x})=0\}.
\end{equation*}
\begin{thm}\label{main theorem without smooth weight}
   Assume that $\mathcal{Z}_{\mathbb{C}}$ is a smooth complete intersection. Then, whenever $n\geq 4R+1$,  one has
    \begin{equation}\label{asymptotic formula without smooth weight}
       N(T)=\mathfrak{S}\mathfrak{J}T^{n-2R}+O(T^{n-2R-\delta}),\ \text{as}\ T\rightarrow \infty, 
    \end{equation}
    for some $\delta>0$, where the singular series $\mathfrak{S}$ is defined in \eqref{eq:singularseriesdef} and the singular integral $\mathfrak{J}$ is defined in \eqref{eq:singularintegraldef}. Furthermore, the constant $\mathfrak{S}\mathfrak{J}$ is positive, provided that the system \eqref{systemofquadraticforms} has a non-trivial solution over $\R$ and $\Z_p$ for all primes $p$.
\end{thm}

\bigskip

We also consider a weighted counting function. Let $\omega:\R^n\to[0,1]$ be a smooth, compactly supported approximation of the characteristic function of the unit ball  $\{\mathbf{x}\in \R^n:\|\mathbf{x}\|\leq 1\}$, satisfying $\omega\equiv1$ in a neighborhood of the origin. For $T>0$ we define the counting function
\begin{equation}\label{eq:smoothcountingdef}
    N_{\omega}(T):=\sum_{\mathbf{x}\in\mathbb{Z}^n}\omega\left(\frac{\mathbf{x}}{T}\right)\prod_{i=1}^R 1_{Q_i(\mathbf{x})=0}(\mathbf{x}),
\end{equation}
where $\|\cdot\|$ stands for the supremum norm. The following result provides the asymptotic formula for $N_{\omega}(T)$ as $T\rightarrow \infty$, provided that $n\geq 4R+1.$

\begin{thm}\label{main theorem with smooth weight}
  Assume that $\mathcal{Z}_{\mathbb{C}}$ is a smooth complete intersection. Then,  whenever $n\geq 4R+1$,   one has
 \begin{equation}\label{asymptotic formula with smooth weight}
       N_{\omega}(T)=\mathfrak{S}\mathfrak{J}_{\omega}T^{n-2R}+O(T^{n-2R-\delta}),\ \text{as}\ T\rightarrow \infty, 
    \end{equation}
    for some $\delta>0$, where the singular series $\mathfrak{S}$ is defined in \eqref{eq:singularseriesdef} and the singular integral $\mathfrak{J}_{\omega}$ is defined in \eqref{def of J_w} depending on $\omega$. Furthermore, the constant $\mathfrak{S}\mathfrak{J}_{\omega}$ is positive, provided that the system \eqref{systemofquadraticforms} has a non-trivial solution over $\R$ and $\Z_p$ for all primes $p$.
\end{thm}

\begin{remark}
The asymptotic formulas in Theorems~\ref{main theorem without smooth weight} and~\ref{main theorem with smooth weight} also hold for systems
$$Q_i(\mathbf{x})=m_i,\qquad 1\le i\le R,$$
where $\mathbf{m}=(m_1,\ldots,m_R)\in\mathbb{Z}^R$ is fixed, under the same smoothness assumption on the homogeneous system $Q_1=\cdots=Q_R=0$.
More precisely, the local densities are given by \eqref{product of padic densities} and \eqref{real density}, with $Q_i$ replaced by $Q_i-m_i$ in both formulas. Indeed, the proofs carry over with minor modifications, since Fourier inversion introduces the additional factor $\operatorname{e}(-\boldsymbol{\xi}\cdot\mathbf{m})$.
\end{remark}

\bigskip

Conjecture \ref{conj:Meyerforsystem} is known for the cases $R=1,2.$ The case $R=1$ is the classical theorem of Meyer \cite{M84}. For $R=2$, Demyanov \cite{De56} established $p$-adic solubility for every prime $p$, whenever $n\geq 9$. The local-global principle was established by Heath-Brown \cite{HB18} for $n\geq 8$. More recently, Morgan and Skorovogatov \cite{MS24}  extended this Hasse principle result to $n\geq 6$, conditional on the finiteness of the Tate–Shafarevich group.

For $R=3$, the results of Birch--Lewis \cite{BL64} and Schuur \cite{S80} show that any three quadratic forms in $n\geq 13$ variables have common nonzero $\mathbb{Z}_p$-zero for any $p\geq 11$. Thus, one concludes from Theorem \ref{thm:Meyerforsystem} that the $p$-adic solubility for $p=2,3,5,7$ implies the existence of a non-trivial integer solution satisfying $(\ref{systemofquadraticforms})$, whenever $n\geq 13$.

We note that the main theorems above reach the square-root cancellation threshold $n\geq 4R+1$ for all $R\geq 1$, as discussed in \cite[p.~367]{BHB17} (see also \cite[p.~207]{RM18}). Furthermore, the asymptotic formula (\ref{asymptotic formula without smooth weight}) agrees with the Manin--Peyre conjecture; see \cite[pp.~364--365]{BHB17}. To place these results in context, we briefly review earlier quantitative results on integer solutions to systems of quadratic forms. For $R\geq 2$, the previously known unconditional asymptotic results for general smooth complete intersections, without additional structure on the quadratic forms, require more than $4R+1$ variables; sharper bounds are known for special families.

For $R=1$, Birch \cite{B62} established the same asymptotic formula as in Theorem~\ref{main theorem without smooth weight} for $n\geq 5$. Heath-Brown \cite{HB96} subsequently developed a $\delta$-symbol method, obtaining the formula in Theorem~\ref{main theorem with smooth weight} in the same range. His work also treats the cases $n=3$ and $n=4$, where logarithmic factors may occur in the main term.

For $R=2$, Browning and Munshi \cite{BM13} obtained a power saving asymptotic for a suitable weighted count on a special family of singular intersections of two quadrics in $n\geq 9$ variables, whose singular locus consists of two conjugate points over $\mathbb{Q}(i)$. A conductor lowering step already appears in their circle method argument and was subsequently developed in Munshi's nested $\delta$- method; see Munshi \cite[Section 1-2]{M15}. For smooth intersections, Munshi \cite{M15} obtained a weighted asymptotic formula for  $n\geq 11$. Under an additional splitting hypothesis on $Q_1$ and $Q_2$, this range was improved to $n\geq 10$  by Heath-Brown and Pierce \cite{HBP17}, who developed a two-dimensional Kloosterman refinement. More recently, Li, Rydin Myerson and Vishe \cite{LRMV24} introduced a two-dimensional $\delta$-symbol method and established the same formula unconditionally for $n\geq 10$. Under the generalized Lindel\"of hypothesis for Dirichlet $L$-functions, their result extends to $n\geq 9$. There has also been a parallel development over function fields. Over $\mathbb{F}_q(t)$ of odd characteristic, Vishe \cite{V23} constructed a two-dimensional Farey dissection and obtained asymptotic formulae for smooth intersections of two quadrics in at least nine variables. Thus, the square-root cancellation threshold is attained for $R=2$ in the function-field setting.

For general $R\geq 1$, the classical work of Birch \cite{B62} gives the asymptotic formula in Theorem~\ref{main theorem without smooth weight} for $n\geq 2R^2+3R$. More generally, Birch's theorem also applies to singular varieties and yields (\ref{asymptotic formula without smooth weight}) whenever $n> 2R(R+1)+\sigma$, where $\sigma$ denotes the dimension of Birch's singular locus. Further refinements of the geometric hypotheses were obtained by Dietmann \cite{D15} and Schindler \cite{S15}. A substantial improvement was subsequently achieved by Rydin Myerson \cite{RM18} through a repulsion argument. In particular, for a smooth complete intersection of codimension $R$, his result requires only $n\geq 9R$ for the asymptotic formula (\ref{asymptotic formula without smooth weight}). 
%In the diagonal case, a sharper bound follows from the work of Cook \cite{C72}. More precisely, if $Q_1,\ldots,Q_R$ are diagonal quadratic forms and the associated projective variety $\mathcal{Z}$ is a smooth complete intersection of codimension $R$, the same asymptotic formula holds for $n\geq 4R+1$. For the homogeneous generalized quadratic forms considered in \cite[Section 1.1]{BPS24}, our results improve the range for the Hasse principle from $m\geq 9$ to $m\geq 5$, where $m$ is the number of variables over the number field, provided that the descended system defines a smooth complete intersection of the expected codimension.
For diagonal quadratic forms, Cook \cite{C72} already attained the square-root cancellation threshold $n\geq 4R+1$, establishing the same asymptotic formula under the smooth complete intersection assumption. Our results therefore extend this sharp variable bound from diagonal systems to general systems of quadratic forms.

As an application, our results also improve the variable bound for the Hasse principle for the homogeneous generalized quadratic forms considered in \cite[Section 1.1]{BPS24}. More precisely, the bound $m\geq 9$, obtained from
Rydin Myerson's result \cite{RM18}, can be improved to $m\geq 5$, where $m$ is the number of variables over the number field, provided that the descended system defines a smooth complete intersection of the expected codimension. Indeed, the descended system consists of $d$ quadratic forms in $dm$ variables, where $d$ is the degree of the number field over $\mathbb{Q}$.

\subsection{Method and relation to previous work}
Our approach combines the Hardy--Littlewood circle method with techniques from homogeneous dynamics. The key connection is that the quadratic exponential sums arising in the circle method can be viewed as Siegel-theta series. The use of theta functions in the study of values of quadratic forms has a long history, including the works of Bentkus--G\"otze \cite{BG99}, G\"otze \cite{G04}, Marklof \cite{M02,M03}, Buterus--G\"otze--Hille--Margulis \cite{BGHM22}, and Kim--Marklof--Welsh \cite{KMW26}. While these works primarily concern Diophantine inequalities for individual quadratic forms, the present paper applies theta-function methods to systems of quadratic equations. In our approach, we use the bounds for Siegel-theta series established by Marklof--Welsh \cite{MW23} to control the exponential sums in terms of height functions on the space of symplectic lattices.

The required moment estimates for these height functions are then obtained using techniques from homogeneous dynamics. In particular, we follow the method of Eskin--Margulis--Mozes \cite{EMM98}, based on mean-value inequalities for lattice height functions (known as Margulis inequalities), as also employed in \cite{BGHM22}. A crucial feature of this approach is that it yields moment estimates for all exponents $\lambda<\frac{2R}{n}$, leading to square-root cancellation bounds. This improves upon the bounds obtainable by conventional exponential sum methods, including the repulsion method of Rydin Myerson \cite{RM18}, and is essential for reaching the sharp variable bound $n\geq 4R+1$.

Applying these techniques from homogeneous dynamics to systems of quadratic equations requires several additional steps. Unlike the Diophantine inequalities considered in \cite{BG99,G04,BGHM22,KMW26}, for which the main contribution in the Fourier analysis comes from a single neighborhood of the origin, our problem concerns Diophantine equations, where the major arcs consist of neighborhoods of rational points with varying denominators. Consequently, the moment estimates are needed not only to control the minor arcs but also to establish the convergence of the singular series and singular integral appearing in the major arcs analysis.

The main technical difficulty lies in establishing the required moment estimates when $R>1$. In this case, the unipotent action on exterior powers gives rise to higher-degree polynomial maps in several variables. The arguments in \cite{EMM98,BGHM22} are based on estimates for linear polynomials, with \cite{BGHM22} reducing the higher-degree case to linear factors. Such a reduction is no longer available in general in our setting. To address this difficulty, we adapt the strategy recently developed in \cite{KO26}, deriving the required sublevel estimates from algebraic dimension bounds obtained using the smoothness assumption \eqref{eq:geometric-hypothesis}.

\subsection{Organization of the paper}
This paper is organized as follows.

Section~\ref{sec:prelim} introduces Siegel-theta series and recalls their properties and connections to the geometry of numbers, following \cite{MW23}. We then state our main moment estimate (Proposition~\ref{prop:mainmomentestimate}) and outline the proofs of the main theorems. In Section~\ref{sec:circlemethod}, we deduce the main theorems from Proposition~\ref{prop:mainmomentestimate}, using the Hardy-Littlewood circle method. The remaining two sections are devoted to the proof of Proposition~\ref{prop:mainmomentestimate}. In Section~\ref{section:mainmomentestimate} we establish the moment estimate assuming the uniform mean-value estimate in Proposition~\ref{prop:localcontraction}, which is proved in Section~\ref{sec:sublevelestimate}.

\section*{Acknowledgments}
W.K. thanks the University of California, Davis for their hospitality and acknowledges support from KIAS individual grant HP101301. K.Y. gratefully acknowledges support from the University of California, Davis, through the KAP allocation. The authors are grateful to Tim Browning, Junxian Li, Hee Oh, Simon Rydin Myerson and Trevor Wooley for their valuable comments. In particular, the authors thank Trevor Wooley for providing a sharpness example in Example~\ref{example:sharpness}, and Tim Browning for pointing out an application of the main results to generalized quadratic forms.

\section*{AI statement}
All mathematical arguments in this paper, with the exception of Example~\ref{example:smoothness} in Section~\ref{section 6} generated by ChatGPT 5.6 Sol and rewritten by the authors, were developed solely by the authors from March 2025 to August 2026. We also used AI tools to assist with proof checking and polishing the manuscript.

\bigskip

\section{Preliminaries and outline of the proof}\label{sec:prelim}

\bigskip

\subsection{Fourier transform}

Let $\omega:\R^n\to[0,1]$ be as in the introduction: it is smooth and compactly supported, takes values in $[0,1]$, approximates the characteristic function of the unit ball $\{\mathbf{x}\in \R^n:\|\mathbf{x}\|\leq 1\}$, and satisfies $\omega\equiv1$ in a neighborhood of the origin. Recall from \eqref{eq:smoothcountingdef} the definition of the counting function $N_{\omega}(T)$:
$$N_{\omega}(T):=\sum_{\mathbf{x}\in\mathbb{Z}^n}\omega\left(\frac{\mathbf{x}}{T}\right)\prod_{i=1}^R 1_{Q_i(\mathbf{x})=0}(\mathbf{x}).$$
By Fourier orthogonality, we may write
\begin{align*}
    N_{\omega}(T)&=\sum_{\bx\in\Z^n}\omega\left(\frac{\mathbf{x}}{T}\right)\int_{[0,1]^R}\operatorname{e}\left(\sum_{i=1}^R \xi_iQ_i(\bx)\right)d\bxi\\&=\int_{[0,1]^R}\left(\sum_{\bx\in\Z^n}\omega\left(\frac{\mathbf{x}}{T}\right)\operatorname{e}\left(\sum_{i=1}^R \xi_iQ_i(\bx)\right) \right)d\bxi,
\end{align*}
where $\bxi=(\xi_1,\ldots,\xi_R)\in [0,1]^R$. Here and hereafter, $\operatorname{e}(x)$ stands for $e^{2\pi i x}$.

\subsection{Siegel-theta series}

We now rewrite the exponential sum appearing after Fourier inversion as a Siegel-theta series. This subsection provides a dictionary between the circle-method notation and the homogeneous-dynamical notation used later.
For $1\leq i\leq R$, let $A_i$ denote the $n\times n$ symmetric matrix associated to $Q_i$, i.e. $Q_i(\mathbf{x})=\tfrac{1}{2}\bx^{\intercal}A_{i}\bx$. For \(\bxi=(\xi_1,\ldots,\xi_R)\in\R^R\), put
\begin{equation}\label{eq:quadraticformmatrixdefinition}
    A_{\bxi}:=\sum_{i=1}^R \xi_iA_i,
\qquad
Q_{\bxi}(\bx):=\sum_{i=1}^R \xi_iQ_i(\bx)
       =\tfrac12\bx^{\intercal}A_{\bxi}\bx .
\end{equation}
Define 
\begin{equation}\label{eq:exponentialsumdefinition}
    S_{\omega}(\bxi;T)
:=\sum_{\bx\in\Z^n}\omega\left(\frac{\bx}{T}\right)
  \operatorname e(Q_{\bxi}(\bx)).
\end{equation}
This is precisely the weighted Siegel-theta series
\(\theta_\omega(T,A_{\bxi},0,0)\) in the notation of Marklof--Welsh \cite{MW23}. We also write $S(\bxi;T):=S_{\chi}(\bxi;T)$, where $\chi$ stands for the characteristic function on the box $[-1,1]^n$.  Equivalently,
\begin{equation}\label{eq:exponentialsumdefinition'}
S(\bxi;T)=\sum_{\bx\in \Z^n\cap [-T,T]^n} \operatorname{e}(Q_{\bxi}(\bx)).
\end{equation}
In particular, $N(T)=N_{\chi}(T)$.

To estimate these exponential sums using Siegel-theta series, we work with the space of symplectic lattices. Let
$$G=\operatorname{Sp}(2n,\R),\qquad
\Gamma=\operatorname{Sp}(2n,\Z),
\qquad
X=G/\Gamma .$$
Here $\operatorname{Sp}(2n,\R)$ acts on $\R^{2n}=\R^n\oplus\R^n$ with
respect to the standard symplectic form. Thus a point of $X$ may be viewed
as a unimodular symplectic lattice in $\R^{2n}$. Let
\begin{equation}\label{eq:subgroupsdefinition}
    a(t)=
\begin{pmatrix}
    e^t\operatorname{Id}_n & 0\\
    0 & e^{-t}\operatorname{Id}_n
\end{pmatrix},
\qquad
u(\bxi)=
\begin{pmatrix}
    \operatorname{Id}_n & A_{\bxi}\\
    0 & \operatorname{Id}_n
\end{pmatrix}.
\end{equation}
The matrix \(u(\bxi)\) is symplectic because \(A_{\bxi}\) is symmetric. We set
\begin{equation}\label{eq:symplecticlatticedef}
    \Lambda_{T,\bxi}:=a(\log T)u(\bxi)\Z^{2n}\in X .
\end{equation}
Equivalently, if \((\bp,\bq)\in\Z^n\oplus\Z^n\), then
\begin{equation}\label{eq:lattice-explicit-form}
 a(\log T)u(\bxi)(\bp,\bq)
 =
 \bigl(T(\bp+A_{\bxi}\bq),\,T^{-1}\bq\bigr).
\end{equation}

The passage from $S_{\omega}(\bxi;T)$ to $\Lambda_{T,\bxi}$ relies on the automorphy of the Siegel theta function arising from the Segal--Shale--Weil representation. The theta function is defined on the Jacobi group $\mathbb H_n\rtimes G$, where $\mathbb H_n$ is the Heisenberg group, and satisfies an automorphy relation; see \cite[Theorem~4.1]{MW23}. This automorphy, together with the estimates in \cite[Section~5]{MW23}, leads to a uniform bound for the theta function in terms of a height function on the space of symplectic lattices. We briefly recall the relevant definitions and estimates.

Following \cite[Section~4]{MW23}, we consider the theta function on the Jacobi group, defined by
\begin{equation}\label{eq:grouptheoreticthetadef}
    \Theta_\omega(h,g):=\sum_{\mathbf m\in\mathbb Z^n}\bigl(W(h)R(g)\omega\bigr)(\mathbf m),
\end{equation}
where $W$ is the Schrödinger representation of the Heisenberg group $\mathbb H_n$, and $R$ is the Segal--Shale--Weil representation of $G=\operatorname{Sp}(2n,\mathbb R)$. For $T>0$, $A\in \operatorname{Sym}_n(\mathbb{R})$, and $\mathbf{x},\mathbf{y}\in\mathbb{R}^n$, we also define
\begin{equation}
    \theta_{\omega}(T,A,\mathbf{x},\mathbf{y}):=\sum_{\mathbf{m}\in\mathbb{Z}^n}\omega(T^{-1}(\mathbf{m}+\mathbf{x}))\operatorname{e}\big(\tfrac{1}{2}\bm^{\intercal}A\bm+\bm^{\intercal}\mathbf{y}\big).
\end{equation}
In particular, $S_{\omega}(\bxi;T)=\theta_\omega(T,A_{\bxi},0,0)$. By \cite[(5.29)]{MW23}, we have
\begin{equation}\label{eq:thetafunctionscomparison}
\theta_\omega(T,A_{\bxi},\mathbf{x},\mathbf{y})=T^{\frac{n}{2}}\Theta_\omega\big((\mathbf{x},\mathbf{y},0),u(\bxi)a(-\log T)\big).
\end{equation}

For our purposes, we only need the following uniform bound from \cite[Section~5.1]{MW23}: there is a height function $D$ on the space of symplectic lattices $X$ such that
\begin{equation}\label{eq:Lipschitzprinciple0}
    |\Theta_\omega(h,g)|\ll_{\omega} D(g^{\operatorname{tr}}\mathbb{Z}^{2n})^{\frac{1}{4}}
\end{equation}
for all $g\in G$ and $h\in \mathbb{H}_n$, where we identify the row lattice $\mathbb Z^{2n}g$ with the column lattice $g^{\operatorname{tr}}\mathbb Z^{2n}$. Now let $g=u(\bxi)a(-\log T)$. The corresponding column lattice is orthogonally equivalent to $\Lambda_{T,\bxi}=a(\log T)u(\bxi)\mathbb Z^{2n}$. Consequently, by \eqref{eq:thetafunctionscomparison}, we get
\begin{equation}\label{eq:exponentialsumuniformheightbound}
    |S_{\omega}(\bxi;T)|=|\theta_\omega(T,A_{\bxi},0,0)|\ll_{\omega} T^{\frac n2}D(\Lambda_{T,\bxi})^{\frac{1}{4}}.
\end{equation}
Here, the factor $T^{\frac{n}{2}}$ reflects the natural square-root size of an $n$-dimensional oscillatory sum. The remaining factor $D(\Lambda_{T,\bxi})^{\frac{1}{4}}$ measures the depth of the associated symplectic lattice in the cusp of $X$. In the next subsection, we relate $D$ to the
height functions arising from the geometry of numbers.

\subsection{Geometry of numbers and height functions ($\alpha$-functions)}

We recall the height functions used to measure the depth of a lattice in the
cusp. Let $\Lambda\subset\R^{2n}$ be a unimodular lattice. We say that a subspace $V\subset\R^{2n}$ is $\Lambda$-rational if $V\cap\Lambda$ is a lattice in $V$. For a $\Lambda$-rational subspace $V$, we denote
$$d_\Lambda(V):=\operatorname{covol}(V\cap\Lambda),$$
where the covolume is computed with respect to the Euclidean measure on $V$. For $1\le k\le2n-1$, define
\begin{equation}\label{eq:alpha-k-definition}
\alpha_k(\Lambda):=
\sup_{\substack{V\text{ is }\Lambda\text{-rational}\\ \dim V=k}}
 d_\Lambda(V)^{-1}.
\end{equation}
Thus $\alpha_k(\Lambda)$ is large when $\Lambda$ contains a
$k$-dimensional rational subspace of small covolume. In particular, $\alpha_1(\Lambda)$ is the reciprocal of the length of the shortest nonzero vector, and $\alpha_n(\Lambda)$ detects $n$-dimensional rational subspaces of small covolume. In terms of successive minima $\lambda_1(\Lambda),\ldots, \lambda_{2n}(\Lambda)$, Minkowski's second theorem gives the comparison
\begin{equation}\label{eq:alpha-successive-minima}
\alpha_k(\Lambda)\asymp_n
\left(\prod_{j=1}^k\lambda_j(\Lambda)\right)^{-1} ,
\end{equation}
where the implied constant depends only on $n$.

For a symplectic lattice $\Lambda=g\Z^{2n}$ with $g\in G$, its Euclidean dual lattice $\Lambda^{*}:=(g^{-1})^{\operatorname{tr}}\Z^{2n}$ satisfies $\Lambda^{*}=J\Lambda$, where $J:=\begin{pmatrix}0&\operatorname{Id}_n\\-\operatorname{Id}_n&0\end{pmatrix}$. Since $J$ is orthogonal, the duality between $k$-dimensional and $(2n-k)$-dimensional covolumes implies that
\begin{equation}\label{eq:symplecticsymmetry}
    \alpha_k(\Lambda)= \alpha_{2n-k}(\Lambda)
\end{equation}
for $0\leq k\leq 2n$. Here and hereafter, we use the convention $\alpha_0(\Lambda)=\alpha_{2n}(\Lambda)=1$. Together with \eqref{eq:alpha-successive-minima}, this also gives the coarse monotonicity
\begin{equation}\label{eq:alpha-monotonicity}
\alpha_1(\Lambda)\ll_n\alpha_2(\Lambda)\ll_n\cdots\ll_n\alpha_n(\Lambda),
\qquad
\alpha_{2n-1}(\Lambda)\ll_n\cdots\ll_n\alpha_n(\Lambda).
\end{equation}
Thus every $\alpha_k$ is controlled by $\alpha_n$ up to a constant. This reduction will be useful later, when we prove moment estimates only for exterior degrees up to $n$.

The height $D$ appearing in \eqref{eq:exponentialsumuniformheightbound} is comparable to $\alpha_n^2$ on symplectic lattices, i.e.
\begin{equation}\label{eq:heightfunctionscomparison}
    D(\Lambda)\asymp_n\alpha_n(\Lambda)^{2}
\end{equation}
for any $\Lambda\in X$.
Thus, the only geometric
quantity from the theta estimate in \cite{MW23} that we need is the smallest
covolume among the $n$-dimensional $\Lambda_{T,\bxi}$-rational subspaces.
Combining this with \eqref{eq:exponentialsumuniformheightbound} gives the following.

\begin{lem}\label{lemma2.1}
For any $T\geq 3$ and \(\bxi\in\R^R\), one has
\begin{equation}\label{exponential sum estimate with smooth weight2}
|S_{\omega}(\bxi;T)|\ll_{\omega} T^{\frac n2}\alpha_n(\Lambda_{T,\bxi})^{\frac12}.
\end{equation}
\end{lem}

\begin{proof}
As above, \(S_{\omega}(\bxi;T)=\theta_\omega(T,A_{\bxi},0,0)\) in the notation of
\cite{MW23}. Applying \eqref{eq:exponentialsumuniformheightbound} and the comparison
\(D\asymp_n\alpha_n^{2}\) gives the claim.
\end{proof}

A similar bound holds for $S(\bxi;T)$ up to logarithmic loss.

\begin{lem}\label{lemma2.2}
For any $T\geq 3$ and $\bxi\in\R^R$, one has
\begin{equation}\label{exponential sum estimate with smooth weight}
|S(\bxi;T)|\ll T^{\frac n2}(\log T)^n\alpha_n(\Lambda_{T,\bxi})^{\frac12}.
\end{equation}
\end{lem}

\begin{proof}
    We follow the arguments in \cite[Subsection 5.3]{MW23}. By \cite[Lemma 5.4]{MW23}, we can find a smooth, compactly supported function $\omega_1:(0,2)\to \mathbb{R}_{\geq 0}$ such that
    $$\chi_{(-1,1)}(x)=\sum_{j\geq 0}\big(\omega_1\big(2^j(1+x)\big)+\omega_1\big(2^j(1-x)\big)\big).$$
    Let $\omega_n(x_1,\ldots,x_n):=\prod_{i=1}^{n}\omega_1(x_i)$.
    For a subset $S\subseteq \{1,\ldots,n\}$ and $\mathbf{j}=(j_1,\ldots,j_n)\in\Z_{\geq 0}^n$ we define     $$g_{\mathbf{j},S}=\begin{pmatrix}
        A_{\mathbf{j}}E_S & 0 \\ 0 & A_{\mathbf{j}}^{-1}E_S
    \end{pmatrix}, \qquad \hat{g}_{\mathbf{j},S}=A_{\mathbf{j}}E_S,$$ where $E_S$ is diagonal with $(i,i)$-entry $-1$ if $i\in S$, $+1$ if $i\notin S$, and $A_{\mathbf{j}}=\operatorname{diag}(2^{j_1},\ldots,2^{j_n})$. Then $$\chi_{(-1,1)^n}(\mathbf{x})=\sum_{\mathbf{j}\geq 0}\sum_{S}\omega_n\big(A_{\mathbf{j}}(E_S\mathbf{x}+\mathbf 1)\big).$$
    Let $\omega_{n,\mathbf{j},S}(\mathbf x):=\omega_n\big(A_{\mathbf{j}}(E_S\mathbf{x}+\mathbf 1)\big)$. 
    
    Without loss of generality, we may assume $T\in \mathbb{N}+\frac{1}{2}$. We decompose $S(\bxi;T)$ as follows:
    \begin{equation}\label{eq:boxexponentialdecomp}
        S(\bxi;T)=\sum_{\mathbf{j}\geq 0}\sum_{S\subseteq \{1,\ldots,n\}} S_{\omega_{n,\mathbf{j},S}}(\bxi;T).
    \end{equation}
    For each $\mathbf{j}$ and $S$, by the transformation formula
    \cite[(5.48), (5.50)]{MW23}, we have
    $$\begin{aligned}
    S_{\omega_{n,\mathbf j,S}}(\bxi;T)&=
    T^{\frac{n}{2}}\Theta_{\omega_{n,\mathbf j,S}}\bigl((0,0,0),u(\bxi)a(-\log T)\bigr)\\&=2^{-\frac12(j_1+\cdots+j_n)}T^{\frac{n}{2}}\Theta_{\omega_n}\bigl((T E_S\mathbf{1},-TA_{\bxi}E_S\mathbf{1},0),u(\bxi)a(-\log T)g_{\mathbf{j},S}\bigr).\end{aligned}$$
    Applying \eqref{eq:Lipschitzprinciple0} and \eqref{eq:heightfunctionscomparison}, we get
    \begin{equation}\label{eq:heightcompensation}
            |S_{\omega_{n,\mathbf{j},S}}(\bxi;T)|\ll  2^{-\frac{1}{2}(j_1+\cdots+j_n)}T^{\frac{n}{2}}\alpha_n(g_{\mathbf j,S}^{-1}
    \Lambda_{T,\bxi})^{\frac{1}{2}}\leq T^{\frac{n}{2}}\alpha_n(\Lambda_{T,\bxi})^{\frac{1}{2}}.
    \end{equation} 
    Here the last inequality follows from $\alpha_n(g_{\mathbf j,S}^{-1}\Lambda) \leq 2^{j_1+\cdots+j_n}\alpha_n(\Lambda)$. 
    
    Moreover, if $\displaystyle\max_{1\leq i\leq n}j_i > \log_2 (4T)$ then we have
    \begin{equation}\label{eq:largeindexvanishing}
        S_{\omega_{n,\mathbf{j},S}}(\bxi;T)=\sum_{\mathbf{x}\in \mathbb{Z}^n}\omega_n\left(A_{\mathbf{j}}\left(\frac{E_S\mathbf{x}}{T}+\mathbf 1\right)\right)\operatorname{e}(Q_{\mathbf{\bxi}}(\mathbf x))=0.
    \end{equation}
    Combining \eqref{eq:boxexponentialdecomp}, \eqref{eq:heightcompensation}, and \eqref{eq:largeindexvanishing}, we conclude that
    \begin{equation*}
        \begin{aligned}
            |S(\bxi;T)|&\leq \sum_{\max j_i\leq \log_2 (4T)}\sum_{S\subseteq \{1,\ldots,n\}}|S_{\omega_{n,\mathbf{j},S}}(\bxi;T)|\\&\ll \sum_{ \max j_i\leq \log_2 (4T)}\sum_{S\subseteq \{1,\ldots,n\}} T^{\frac{n}{2}}\alpha_n(\Lambda_{T,\bxi})^{\frac{1}{2}}\ll T^{\frac{n}{2}}(\log T)^n \alpha_n(\Lambda_{T,\bxi})^{\frac{1}{2}}.
        \end{aligned}
    \end{equation*}
\end{proof}

\begin{remark}
The following variant of Lemma~\ref{lemma2.2} holds uniformly for translated boxes:
\begin{equation}\label{eq:lem2.2variant}
    \left|\sum_{\bx\in \Z^n\cap ([-T,T]^n+\mathbf{x}_0)} \operatorname{e}(Q_{\bxi}(\bx))\right|\ll T^{\frac n2}(\log T)^n\alpha_n(\Lambda_{T,\bxi})^{\frac12}
\end{equation}
    for any $\mathbf{x}_0\in\mathbb{R}^n$. This follows by the same argument as in the proof of Lemma~\ref{lemma2.2}, using \eqref{eq:thetafunctionscomparison} with a general translation parameter $\mathbf{x}$. Here, the endpoints of the translated box may be adjusted to half-integers without changing its lattice points.
\end{remark}

The following moment bound is the main
ingredient of our proof:
\begin{prop}\label{prop:mainmomentestimate}
Assume that $\mathcal{Z}_\C$ is a smooth complete intersection and that
$n\geq4R+1$.
For any $0\leq \lambda < \frac{2R}{n}$ and any $\varepsilon>0$, there exists a
constant $0<C_{\lambda,\varepsilon}<\infty$ such that
\begin{equation}\label{moment bound}
    \int_{[0,1]^R}\alpha_n(\Lambda_{T,\bxi})^{\lambda}d\bxi
    \leq C_{\lambda,\varepsilon}T^{\varepsilon}
\end{equation}
for all $T\geq 1$.
\end{prop}

\subsection{Outline of the proof}
We apply the Hardy--Littlewood circle method to the identity
\begin{equation}\label{eq:circlemethodexpression}
    N(T)=\int_{\mathbb{T}^R}S(\bxi;T)\,d\bxi=\int_{\mathfrak{M}_{T,\delta}}S(\bxi;T)\,d\bxi+\int_{\mathfrak{m}_{T,\delta}}S(\bxi;T)\,d\bxi,
\end{equation}
decomposing the frequency space $\mathbb{T}^R$ into the major arcs $\mathfrak{M}_{T,\delta}$ and the minor arcs $\mathfrak{m}_{T,\delta}=[0,1]^R\setminus \mathfrak{M}_{T,\delta}$. The major arcs are neighborhoods of rational points with small denominators and their contribution will give the main term $c T^{n-2R}$, where $c\geq 0$ is the product of local densities. Here, the exponent $n-2R$ reflects quadratic homogeneity. Indeed, near $\bxi=0$ the exponential sum has size $T^n$ on the natural frequency scale $T^{-2}$ in each of $R$ coordinates. Integrating over this region therefore yields the scale $T^n(T^{-2})^R=T^{n-2R}$. Summing the corresponding contributions near other rational points accounts for the congruence conditions and produces the $p$-adic density factors in $c$.

The natural square-root cancellation scale for an exponential sum with $\asymp T^n$ terms is $T^{\frac{n}{2}}$. Comparing this with the main-term scale $T^{n-2R}$ suggests the condition $n-2R>\frac{n}{2}$, or equivalently $n>4R$. To establish an asymptotic formula for $N(T)$ in view of \eqref{eq:circlemethodexpression}, we must control the contribution of large values of $S(\bxi;T)$. Lemma~\ref{lemma2.2} bounds this exponential sum in terms of a lattice height:
$$|S(\boldsymbol{\xi};T)|
\ll T^{\frac{n}{2}}(\log T)^n\alpha_n(\Lambda_{T,\boldsymbol\xi})^{\frac{1}{2}},$$
where $\Lambda_{T,\boldsymbol\xi}=a(\log T)u(\boldsymbol\xi)\mathbb Z^{2n}$. Thus the exponential sum is controlled by the geometry of the associated symplectic lattice $\Lambda_{T,\boldsymbol\xi}$.

Proposition~\ref{prop:mainmomentestimate} bounds the moments of $\alpha_n(\Lambda_{T,\boldsymbol\xi})$ over $\bxi$. This estimate enters both the major and minor arcs analysis. In \eqref{eq:circlemethodexpression}, the contribution of the major arcs is initially approximated by the product of a truncated singular series and a truncated singular integral. Proposition~\ref{prop:mainmomentestimate} controls the tails of both factors. This allows us to replace the truncated factors by their limits with a power-saving error, yielding the main term $c T^{n-2R}$, where $c$ is the product of local densities. On the minor arcs, combining Proposition~\ref{prop:mainmomentestimate} with a pointwise saving (Lemma~\ref{pointwise bounds}) gives a contribution of order $O(T^{n-2R-\tau})$ for some $\tau>0$.

We now turn to the proof of Proposition~\ref{prop:mainmomentestimate}. We follow the method of Eskin–Margulis–Mozes \cite{EMM98}, which proceeds in three steps. First, we establish a uniform mean-value estimate for individual vectors in exterior powers of $\mathbb{R}^{2n}$, stated in Proposition~\ref{prop:localcontraction}. We then deduce a Margulis inequality for a suitable weighted sum of the lattice height functions $\alpha_j^\lambda$. Finally, we iterate this inequality at a sufficiently large fixed averaging time to obtain the required $T^\varepsilon$ bound. The latter two steps follow the argument of \cite{EMM98} with minor adaptations and are carried out in Section~\ref{section:mainmomentestimate}.

The main additional difficulty lies in the first step, namely establishing Proposition~\ref{prop:localcontraction} for $R>1$. To obtain the uniform mean-value estimate in the required range $0<\lambda<\frac{2R}{n}$, we need suitable integrability estimates for the polynomial maps arising from the unipotent action on exterior powers. Following the strategy of \cite{KO26}, we derive these estimates from algebraic dimension bounds. More precisely, in Section~\ref{sec:sublevelestimate}, we first use the smoothness of $\mathcal Z_{\mathbb C}$ to establish the required dimension estimates, and then deduce uniform polynomial sublevel estimates that imply Proposition~\ref{prop:localcontraction}.

\bigskip

\section{The circle method}\label{sec:circlemethod}

We split the domain $[0,1]^R$ into two regions as follows. Let $0<\delta<1$ and set
$$\mathfrak{M}_{T,\delta}:=\bigcup_{\substack{q\in\mathbb{N}\\ q\leq T^{\delta}}}\bigcup_{\substack{0\leq p_1,\ldots,p_R<q\\ (p_1,\ldots,p_R,q)=1}}\mathfrak{M}_{T,\delta}(q,\mathbf{p}),$$
where $$\mathfrak{M}_{T,\delta}(q,\mathbf{p}):=\left\{\bxi\in[0,1]^R:\left\|\xi_i-\frac{p_i}{q}\right\|_{\mathbb{Z}}<T^{\delta-2},\ 1\leq i\leq R\right\}.$$
Here $\|s\|_{\mathbb Z}:=\min_{a\in\Z}|s-a|$ for $s\in\R$. Furthermore, we define
$$\mathfrak{m}_{T,\delta}=[0,1]^R\setminus \mathfrak{M}_{T,\delta}.$$
Throughout this section, differences between frequency parameters are taken modulo $\mathbb{Z}^R$, with representatives chosen in the corresponding small neighborhoods of the origin.

We record the following elementary pointwise upper bounds:
\begin{lem}\label{pointwise bounds}
    Let $T>1$. We have
    \begin{enumerate}
        \item $\alpha_n(\Lambda_{T,\bxi})\leq T^n$ for any $\bxi\in[0,1]^R$,
        \item $\alpha_n(\Lambda_{T,\bxi})\ll T^{n-\delta}$ for any $\bxi\in \mathfrak{m}_{T,\delta}$.
    \end{enumerate}
\end{lem}

For~(2), compare the analogous estimate for exponential sums in \cite[Lemma~2.3]{RM18}.

\begin{proof}
Write $\Lambda=\Lambda_{T,\bxi}$. Let $W$ be a
$\Lambda$-rational subspace with $\dim W=n$. By
\eqref{eq:lattice-explicit-form}, a $\Z$-basis of $W\cap\Lambda$
forms the columns of a matrix
\begin{equation}\label{eq:Mdefinition}
M=\begin{pmatrix}T(Y+A_{\bxi}Z)\\T^{-1}Z\end{pmatrix},    
\end{equation}
where $Y,Z\in \operatorname{Mat}_{n\times n}(\Z)$ and the stacked
matrix $\left(\begin{smallmatrix}Y\\Z\end{smallmatrix}\right)$
has rank $n$. By the Cauchy--Binet formula, every $n\times n$
minor of $M$ has absolute value at most $d_\Lambda(W)$.

Let $r=\operatorname{rank}Z$. Choose $r$ linearly independent rows
of $Z$ and complete them to a basis of $\R^n$ using $n-r$ rows
of $Y$. These rows form an integer matrix with nonzero determinant.
Every row of $A_{\bxi}Z$ is a linear combination of the selected
rows of $Z$, so the $A_{\bxi}Z$ terms do not change this determinant.
The corresponding minor of $M$ is therefore $T^{n-2r}$ times
a nonzero integer. Hence
\[
d_\Lambda(W)\ge T^{n-2r}\ge T^{-n}.
\]
Taking the supremum of $d_\Lambda(W)^{-1}$ over $W$ proves~(1).

For~(2), assume now that
$\bxi\in\mathfrak m_{T,\delta}$. If $r<n$, the bound
$d_\Lambda(W)\ge T^{n-2r}\geq T^{2-n}$ already gives the required estimate,
since $\delta<1$. It therefore remains to consider $r=n$. Let $q_0:=|\det Z|\in\N$ and $Y':=-q_0YZ^{-1}\in\operatorname{Mat}_{n\times n}(\Z)$. The bound for the bottom minor gives \(q_0\le d_\Lambda(W)T^n\). Factoring $Z$ on the right of \eqref{eq:Mdefinition}, we can write
\begin{equation}\label{eq:Mfactoring}
    M=\begin{pmatrix}
T(A_{\bxi}+YZ^{-1})\\
T^{-1}I_n
\end{pmatrix}Z.
\end{equation}
For $1\le i,j\le n$, the minor consisting of the $i$th top row
and all bottom rows except the $j$th therefore equals, up to sign,
$T^{2-n}(q_0A_{\bxi}-Y')_{ij}$. Consequently,
\begin{equation}\label{eq:Y'approximation}
    \max_{i,j}|(q_0A_{\bxi}-Y')_{ij}|\le d_\Lambda(W)T^{n-2}.
\end{equation}

The smoothness assumption \eqref{eq:geometric-hypothesis} implies that
$A_1,\ldots,A_R$ are linearly independent over $\Q$.
Thus, given $A_1,\ldots,A_R$, there exist $m\in \mathbb{N}$ and linear forms
$\psi_1,\ldots,\psi_R$ on $\operatorname{Mat}_{n\times n}(\R)$,
with integer coefficients such that $\psi_i(A_j)=m\delta_{ij}$ for $1\leq i,j\leq R$. Let us write $q=mq_0$ and
$$\bp=(\psi_1(Y'),\ldots,\psi_R(Y'))\in\Z^R.$$
Then we have
\begin{equation}\label{eq:largedenominatorcase}
    q\le md_\Lambda(W)T^n,\qquad
\max_i|q\xi_i-p_i|
=\max_i|\psi_i(q_0A_{\bxi}-Y')|
\le md_\Lambda(W)T^{n-2}.
\end{equation}
Note that $m$ depends only on the fixed matrices $A_1,\ldots,A_R$.

If $q\le T^\delta$, then by the definition
of $\mathfrak m_{T,\delta}$, we also have
\begin{equation}\label{eq:smalldenominatorcase}
    \max_i\left\|\xi_i-\frac{p_i}{q}\right\|_{\mathbb Z}
\ge T^{\delta-2}.
\end{equation}
Thus, whether or not $q\le T^\delta$, \eqref{eq:largedenominatorcase} and \eqref{eq:smalldenominatorcase} give
$$1\le\max\left\{
qT^{-\delta},
T^{2-\delta}\max_i
\left\|\xi_i-\frac{p_i}{q}\right\|_{\mathbb Z}
\right\}\le md_\Lambda(W)T^{n-\delta}.$$
Together with the case $r<n$, it follows that
$d_\Lambda(W)^{-1}\ll T^{n-\delta}$ for every $W$. Taking the supremum over $W$ proves~(2).
\end{proof}

\bigskip

\subsection{Major arcs analysis}\label{subsec 3.1}
In this subsection, we shall provide the asymptotic formulas for the mean values of $S(\bxi;T)$ and $S_{\omega}(\bxi;T)$ over major arcs, given by
\begin{equation}\label{mean value of the exponential sum over major arcs}
\begin{aligned}
\int_{\mathfrak{M}_{T,\delta}}S(\bxi;T)d\bxi\ \ \text{and }\int_{\mathfrak{M}_{T,\delta}}S_{\omega}(\bxi;T)d\bxi.
\end{aligned}
\end{equation}
We note that the asymptotic formula for the second mean value can be deduced in a similar manner at the end of Subsection \ref{subsec 3.1}. Therefore, the majority of this subsection is devoted to evaluating $$\int_{\mathfrak{M}_{T,\delta}}S(\bxi;T)d\bxi.$$

\bigskip

Using a variant of the moment estimate (\ref{moment bound}) combined with a routine major arcs estimate, we have the following:
\begin{lem}\label{asymptotic formula for the mean value over major arcs}
    Let $n\geq 4R+1$ and fix $0<\delta<\frac{1}{2R+3}.$ Then one has
    \begin{equation*}
    \int_{\mathfrak{M}_{T,\delta}}S(\bxi;T)d\bxi=\mathfrak{S}\mathfrak{J}T^{n-2R}+O(T^{n-2R-\tau}),
    \end{equation*}
   for some $\tau>0,$ where the singular series $\mathfrak{S}$ is defined in \eqref{eq:singularseriesdef} depending on $Q_1,\ldots,Q_R$ and the singular integral $\mathfrak{J}$ is defined in \eqref{eq:singularintegraldef} depending on  $Q_1,\ldots,Q_R$. Furthermore, the constant $\mathfrak{S}\mathfrak{J}$ is positive, provided that the system \eqref{systemofquadraticforms} has a non-trivial solution over $\R$ and $\Z_p$ for all primes $p$.
\end{lem}

We will prove Lemma $\ref{asymptotic formula for the mean value over major arcs}$ at the end of Subsection $\ref{subsec 3.1}$, using the following lemmas \ref{mean value estimate over major arcs}, \ref{absolute convergence of the singular series} and $\ref{absolute convergence of singular integral}$. We then establish the corresponding weighted version (Lemma $\ref{asymptotic formula of mean value over major arcs 2}$). To describe the local functions $\mathfrak{S}$ and $\mathfrak{J}$, we introduce some definitions.
Define 
\begin{equation*}
    \begin{aligned}
        S_q(\mathbf{p}):=\sum_{\bx\in (\Z/q\Z)^n}e\left(\sum_{i=1}^R\frac{p_i}{q}Q_i(\bx)\right)
         \end{aligned}
\end{equation*}
for $q\geq 1$ and $\mathbf{p}=(p_1,\ldots,p_R)\in\mathbb{Z}^R$ with $1\leq p_1,\ldots,p_R\leq q$, and 
\begin{equation*}
\begin{aligned}I(\boldsymbol{\beta};T):=\int_{|\boldsymbol{\gamma}|\leq T}e\left(\sum_{i=1}^R\beta_iQ_i(\boldsymbol{\gamma})\right)d\boldsymbol{\gamma}
    \end{aligned}
\end{equation*}
for $\boldsymbol{\beta}\in\mathbb{R}^R$ and $T\geq 1$. With these functions, we define the local functions
\begin{equation}\label{eq:singularseriesdef}
    \begin{aligned}
\mathfrak{S}&:=\sum_{q=1}^{\infty}\sum_{\substack{0\leq p_1,\ldots,p_R< q\\(q,\mathbf{p})=1}}q^{-n}S_q(\mathbf{p}) \end{aligned}
\end{equation}
and
\begin{equation}
    \begin{aligned}\label{eq:singularintegraldef}
\mathfrak{J}&:=\int_{\boldsymbol{\beta}\in \mathbb{R}^R}I(\boldsymbol{\beta};1)d\boldsymbol{\beta}
    \end{aligned}
\end{equation}

\bigskip

For the proof of Lemma \ref{asymptotic formula for the mean value over major arcs},  we first deduce the asymptotic formula for the mean value of $S(\boldsymbol{\xi},T)$ over the major arcs in terms of the truncated singular series and integrals. Before stating the result, we define the truncated local functions 
\begin{equation*}
\begin{aligned}
    \mathfrak{S}(B):=\sum_{1\leq q\leq B}\sum_{\substack{0\leq p_1,\ldots,p_R< q\\(q,\mathbf{p})=1}}q^{-n}S_q(\mathbf{p})\end{aligned}
\end{equation*}
and
\begin{equation*}
\begin{aligned}
    \mathfrak{J}(B;T):=\int_{\|\boldsymbol{\beta}\|_{\infty}\leq BT^{-2}}I(\boldsymbol{\beta};T)d\boldsymbol{\beta}
\end{aligned}
\end{equation*}
for $B,T\geq 1$.

\begin{lem}\label{mean value estimate over major arcs}
    When $0<\delta<\frac{1}{2R+3},$ one has
    \begin{equation*}
\int_{\mathfrak{M}_{T,\delta}}S(\bxi;T)d\bxi=\mathfrak{S}(T^{\delta})\mathfrak{J}(T^{\delta};1)T^{n-2R}+O(T^{n-2R-\tau}),
    \end{equation*}
    for some $\tau>0.$
\end{lem}
\begin{proof}
    Let $1\leq q\leq T^{\delta}$ and $0\leq p_1,\ldots,p_R<q$. For $\bxi=(\xi_1,\ldots,\xi_R)\in \mathfrak{M}_{T,\delta}(q,\mathbf{p}),$ write 
    \begin{equation*}
        \boldsymbol{\beta}=(\beta_1,\ldots,\beta_R),\qquad \beta_i=\xi_i-\frac{p_i}{q},\ \ (1\leq i\leq R).
    \end{equation*}
    On writing $\bx=q\mathbf{y}+\mathbf{r}$ with $\mathbf{y}\in\mathbb{Z}^n$ and $\mathbf{r}\in\{0,\ldots,q-1\}^n$ in the definition \eqref{eq:exponentialsumdefinition'} of $S(\bxi;T)$, we see that 
    \begin{equation}\label{splitting S(alpha) into two} 
        \begin{aligned}
S(\bxi;T)&=\sum_{\mathbf{r}\in\{0,\ldots,q-1\}^n}\sum_{\substack{\mathbf{y}\in\mathbb{Z}^n\\|q\mathbf{y}+\mathbf{r}|\leq T}}e\bigg(\sum_{i=1}^R\bigg(\beta_i+\frac{p_i}{q}\bigg)Q_i(q\mathbf{y}+\mathbf{r})\bigg)\\
&=\sum_{\mathbf{r}\in\{0,\ldots,q-1\}^n}e\bigg(\sum_{i=1}^R\frac{p_i}{q}Q_i(\mathbf{r})\bigg)\cdot \sum_{\substack{\mathbf{y}\in\mathbb{Z}^n\\|q\mathbf{y}+\mathbf{r}|\leq T}}e\bigg(\sum_{i=1}^R\beta_iQ_i(q\mathbf{y}+\mathbf{r})\bigg).
        \end{aligned}
    \end{equation}

By the standard argument using the mean value theorem, as in \cite[Lemma 5.1]{B62}, one has
    \begin{equation}\label{the standard mean value theorem}
\sum_{\substack{\mathbf{y}\in\mathbb{Z}^n\\|q\mathbf{y}+\mathbf{r}|\leq T}}e\bigg(\sum_{i=1}^R\beta_iQ_i(q\mathbf{y}+\mathbf{r})\bigg)=q^{-n}I(\boldsymbol{\beta};T)+O\bigg(q^{-n}T^{n-1}\big(q+\sum_{i=1}^RqT^2|\beta_i|\big)\bigg).
    \end{equation}
    Substituting this into the last expression of (\ref{splitting S(alpha) into two}), we see that whenever $\bxi\in \mathfrak{M}_{T,\delta}(q,\mathbf{p}),$ one has
    \begin{equation}\label{approximation in major arcs}
        S(\bxi;T)=q^{-n}S_q(\mathbf{p})I(\boldsymbol{\beta};T)+O(qT^{n-1+\delta}).
    \end{equation}
    Hence, we deduce that 
    \begin{equation}\label{last step of the proof of Lemma 3.3}
        \begin{aligned}
\int_{\mathfrak{M}_{T,\delta}}S(\bxi;T)d\bxi&=\sum_{1\leq q\leq T^{\delta}}\sum_{\substack{0\leq p_1,\ldots,p_R< q \\ (q,\mathbf{p})=1}}\int_{\mathfrak{M}_{T,\delta}(q,\mathbf{p})}S(\bxi;T)d\bxi\\
&=\mathfrak{S}(T^{\delta})\mathfrak{J}(T^{\delta};T)+O(T^{n-2R-1+\delta(2R+3)})\\
&=\mathfrak{S}(T^{\delta})\mathfrak{J}(T^{\delta};1)T^{n-2R}+O(T^{n-2R-1+\delta(2R+3)}).
        \end{aligned}
    \end{equation}
    Since $0<\delta<\frac{1}{2R+3}$, this completes the proof of Lemma $\ref{mean value estimate over major arcs}$.
\end{proof}

\bigskip

Next,  we show that $\mathfrak{S}(B)$ and $\mathfrak{J}(B;1)$ converge absolutely as $B\to\infty$. Furthermore, these converge to the singular series $\mathfrak{S}$ and the singular integral $\mathfrak{J}$, respectively. 
\begin{lem}\label{absolute convergence of the singular series}
 Whenever $n\geq 4R+1$, 
 the truncated singular series $\mathfrak{S}(B)$ converges absolutely, as $B\rightarrow \infty$. Furthermore, one has
 \begin{equation*}
     |\mathfrak{S}-\mathfrak{S}(B)|\ll B^{-\tau} \ \ \text{and}\ \ \mathfrak{S}\ll 1,
 \end{equation*}
 for some $\tau>0.$
\end{lem}
\begin{proof}
By definition of $\mathfrak{S}$ and $\mathfrak{S}(B)$, we may write
\begin{equation}\label{eq:truncationdecomposition}
\begin{aligned}
    |\mathfrak{S}-\mathfrak{S}(B)|&=\big|\sum_{q> B}\sum_{\substack{0\leq p_1,\ldots,p_R<q\\(q,\mathbf{p})=1}}q^{-n}S_q(\mathbf{p})\big|\\&\leq \sum_{i=0}^{\infty}\sum_{2^{i}B< q\leq 2^{i+1}B}\sum_{\substack{0\leq p_1,\ldots,p_R< q\\(q,\mathbf{p})=1}}q^{-n}|S_q(\mathbf{p})|.
\end{aligned}
\end{equation}
For any $T<q\leq 2T$ and $\mathbf{p}$, we deduce from \eqref{eq:lem2.2variant} that
\begin{equation}\label{eq:Sqpalphacomparison}
\begin{aligned}
        |S_q(\mathbf{p})|&=\left|\sum_{\mathbf{x}\in \mathbb{Z}^n \cap [\frac{1}{2},q+\frac{1}{2}]^n} \operatorname{e}(Q_{\mathbf{p}/q}(\mathbf{x}))\right|\\&\ll q^{\frac{n}{2}}(\log q)^n\alpha_n(\Lambda_{q/2,\mathbf{p}/q})^{\frac{1}{2}}\asymp T^{\frac{n}{2}}(\log T)^n \alpha_n(\Lambda_{T,\mathbf{p}/q})^{\frac{1}{2}}.
\end{aligned}
\end{equation}
Here $\alpha_n(\Lambda_{q/2,\mathbf{p}/q})\asymp \alpha_n(\Lambda_{T,\mathbf{p}/q})$ follows from $\Lambda_{T,\mathbf{p}/q}=a\big(\log \big(\frac{2T}{q}\big)\big)\Lambda_{q/2,\mathbf{p}/q}$ and that $\log\big(\frac{2T}{q}\big)$ is bounded. Moreover, for any $\varepsilon>0$ we have $\mathbf{p}/q\in \mathfrak{m}_{T,\frac{1}{2}-\varepsilon}$ for sufficiently large $q$. Hence (2) of Lemma~\ref{pointwise bounds} gives
\begin{equation}\label{majorarcpointwisebound}
    \alpha_n(\Lambda_{T,\mathbf{p}/q})\ll_{\varepsilon} T^{n-\frac{1}{2}+\varepsilon}.
\end{equation}

    Let $0<\lambda<\frac{2R}{n}$ and $\varepsilon>0$. From Proposition \ref{prop:mainmomentestimate}, for any $T>1$ we have
\begin{equation}\label{application of moment estimates in major arcs}
\begin{aligned}
    C_{\lambda,\varepsilon}T^{\varepsilon}&\geq \int_{[0,1]^R}\alpha_n(\Lambda_{T,\bxi})^{\lambda} d \bxi\\
    &\geq \sum_{T< q\leq 2T}\sum_{\substack{0\leq p_1,\ldots,p_R<q\\ (p_1,\ldots,p_R,q)=1}}\int_{B(\frac{\mathbf{p}}{q},\frac{1}{10T^2})}\alpha_n(\Lambda_{T,\bxi})^{\lambda}d \bxi.
\end{aligned}
\end{equation}
On noting that $\alpha_n(\Lambda_{T,\bxi})\asymp \alpha_n(\Lambda_{T,\bxi'})$, provided that $\|\bxi-\bxi'\|\ll T^{-2},$ we see from (\ref{application of moment estimates in major arcs}) that
\begin{equation}\label{discrete moment estimate}
    C_{\lambda,\varepsilon}T^{\varepsilon}\gg \sum_{T< q\leq 2T}\sum_{\substack{0\leq p_1,\ldots,p_R<q\\ (p_1,\ldots,p_R,q)=1}}T^{-2R}\alpha_n(\Lambda_{T,\mathbf{p}/q})^{\lambda},
\end{equation}
where we used the fact that $\text{Leb}\big(B(\frac{\mathbf{p}}{q},\frac{1}{10T^2})\big)\asymp T^{-2R}$ for each $(q,\mathbf{p})$.

Combining \eqref{eq:Sqpalphacomparison}, \eqref{majorarcpointwisebound}, and \eqref{discrete moment estimate}, it follows that
\begin{equation}
\begin{aligned}
    &\sum_{T< q\leq 2T}q^{-n}\sum_{\substack{0\leq p_1,\ldots,p_R<q\\ (p_1,\ldots,p_R,q)=1}}|S_q(\mathbf{p})|\\ &\ll \sum_{T< q\leq 2T}q^{-n}\sum_{\substack{0\leq p_1,\ldots,p_R<q\\ (p_1,\ldots,p_R,q)=1}}T^{n/2}(\log T)^n\alpha_n(\Lambda_{T,\mathbf{p}/q})^{1/2}\\
&\ll T^{-n/2}(\log T)^n\sup_{\substack{T< q\leq 2T\\0\leq p_1,\ldots,p_R<q\\ (p_1,\ldots,p_R,q)=1}}\alpha_n(\Lambda_{T,\mathbf{p}/q})^{\frac{1}{2}-\lambda}\cdot \sum_{T< q\leq 2T}\sum_{\substack{0\leq p_1,\ldots,p_R<q\\ (p_1,\ldots,p_R,q)=1}}\alpha_n(\Lambda_{T,\mathbf{p}/q})^{\lambda}\\
&\ll_{\varepsilon} T^{-n/2}(\log T)^n(T^{n-\frac{1}{2}+\varepsilon})^{\frac{1}{2}-\lambda}\cdot \sum_{T< q\leq 2T}\sum_{\substack{0\leq p_1,\ldots,p_R<q\\ (p_1,\ldots,p_R,q)=1}}\alpha_n(\Lambda_{T,\mathbf{p}/q})^{\lambda}\\
&\ll_{\varepsilon,\lambda} T^{-n/2}(\log T)^n(T^{n-\frac{1}{2}+\varepsilon})^{\frac{1}{2}-\lambda} T^{2R+\varepsilon}\\&= T^{-\frac{1}{2}(\frac{1}{2}-\frac{2R}{n})+(n-\frac{1}{2}+\varepsilon)(\frac{2R}{n}-\lambda)+\varepsilon(\frac{3}{2}-\frac{2R}{n})}(\log T)^n.
\end{aligned}
\end{equation}
Since $\varepsilon>0$ can be taken arbitrarily small and $\lambda<\frac{2R}{n}$ arbitrarily close to $\frac{2R}{n}$, we get
\begin{equation}\label{eq:majorarcatscaleT}
    \sum_{T< q\leq 2T}q^{-n}\sum_{\substack{0\leq p_1,\ldots,p_R<q\\ (p_1,\ldots,p_R,q)=1}}|S_q(\mathbf{p})|\ll_{\tau} T^{-\tau} 
\end{equation}
for any $0<\tau<\frac{1}{2}(\frac{1}{2}-\frac{2R}{n})$. Summing \eqref{eq:majorarcatscaleT} over $T=2^iB$, we conclude that $|\mathfrak{S}-\mathfrak{S}(B)|\ll B^{-\tau}$, in view of \eqref{eq:truncationdecomposition}. Furthermore, the inequality $
    \mathfrak{S}\ll 1$ immediately follows by $(\ref{eq:majorarcatscaleT}).$
\end{proof}

\bigskip

\begin{lem}\label{absolute convergence of singular integral}
    Whenever $n\geq 4R+1$, the truncated singular integral $\mathfrak{J}(B;1)$ converges absolutely, as $B\rightarrow \infty$. Furthermore, one has
 \begin{equation*}
     |\mathfrak{J}-\mathfrak{J}(B;1)|\ll B^{-\tau}\ \ \text{and}\ \ \mathfrak{J}\ll 1,
 \end{equation*}
 for some $\tau>0.$
\end{lem}
\begin{proof}
    By definition of $\mathfrak{J}$ and $\mathfrak{J}(B;1)$ we have
    $$\mathfrak{J}-\mathfrak{J}(B;1)=\sum_{i=0}^{\infty}\int_{2^i B\leq \|\boldsymbol{\beta} \|_{\infty}\leq 2^{i+1}B}I(\boldsymbol{\beta};1)d\boldsymbol{\beta},$$
    hence it suffices to show that
    \begin{equation*}
        \int_{B\leq \|\boldsymbol{\beta}\|_{\infty}\leq 2B}|I(\boldsymbol{\beta};1)|d\boldsymbol{\beta}\ll B^{-\tau},
    \end{equation*}
    for some $\tau>0.$
   
    Define $T$ to be a positive number satisfying $B=T^{\delta}$ with $\delta <\frac{1}{R+1}$. Then, by the approximation in $(\ref{approximation in major arcs})$ with $q=1$, one infers that
    \begin{equation}\label{approximation for the singular integrals}
    \begin{aligned}
S(T^{-2}\boldsymbol{\beta};T)&=I(T^{-2}\boldsymbol{\beta};T)+O(T^{n-1+\delta})\\
&=T^{n}I(\boldsymbol{\beta};1)+O(T^{n-1+\delta}).
    \end{aligned}
    \end{equation}
    Therefore, we find from (\ref{approximation for the singular integrals}) and (\ref{exponential sum estimate with smooth weight}) that
    \begin{equation}\label{key step for absolute convergence of the singular integrals}
    \begin{aligned}
          \int_{B\leq \|\boldsymbol{\beta}\|_{\infty}\leq 2B}|I(\boldsymbol{\beta};1)|d\boldsymbol{\beta}&= \int_{B\leq \|\boldsymbol{\beta}\|_{\infty}\leq 2B}\bigr(T^{-n}|S(T^{-2}\boldsymbol{\beta};T)|+O(T^{-1+\delta})\bigr)d\boldsymbol{\beta}\\
          &=T^{2R-n}\int_{\frac{B}{T^2}\leq \|\boldsymbol{\beta}\|_{\infty}\leq \frac{2B}{T^2}}|S(\boldsymbol{\beta};T)|d\boldsymbol{\beta}+O(T^{-1+\delta}B^{R})\\
          &\ll T^{2R-n}T^{\frac{n}{2}}(\log T)^n\int_{\frac{B}{T^2}\leq \|\boldsymbol{\beta}\|_{\infty}\leq \frac{2B}{T^2}}\alpha_n(\Lambda_{T,\boldsymbol{\beta}})^{\frac12}d\boldsymbol{\beta}+O(T^{-1+\delta}B^{R})
    \end{aligned}
    \end{equation}

    One can easily check that whenever  $\frac{B}{T^2}\leq \|\boldsymbol{\beta}\|_{\infty}\leq \frac{2B}{T^2}$, one has $\boldsymbol{\beta}\notin \mathfrak{M}_{T,\delta-\eta}$ for sufficiently small $\eta>0.$ Hence, by Lemma \ref{pointwise bounds}, we see that
    \begin{equation*}
      \sup_{\frac{B}{T^2}\leq \|\boldsymbol{\beta}\|_{\infty}\leq \frac{2B}{T^2}}  \alpha_n(\Lambda_{T,\boldsymbol{\beta}})\ll T^{n-\eta_1},\ \text{for some}\ \eta_1>0.
    \end{equation*}
    Since $n\geq 4R+1$, we may choose a sufficiently small $\epsilon>0$ so that $\frac{1}{2}-\frac{2R}{n}+\epsilon>0$. Therefore, with this $\epsilon>0,$ it follows by Proposition \ref{prop:mainmomentestimate} that  one has
    \begin{equation}
        \begin{aligned}
            \int_{\frac{B}{T^2}\leq \|\boldsymbol{\beta}\|_{\infty}\leq \frac{2B}{T^2}}\alpha_n(\Lambda_{T,\boldsymbol{\beta}})^{\frac12}d\boldsymbol{\beta}&\ll  \sup_{\frac{B}{T^2}\leq \|\boldsymbol{\beta}\|_{\infty}\leq \frac{2B}{T^2}}  (\alpha_n(\Lambda_{T,\boldsymbol{\beta}}))^{\epsilon+\frac{1}{2}-\frac{2R}{n}} \int_{\frac{B}{T^2}\leq \|\boldsymbol{\beta}\|_{\infty}\leq \frac{2B}{T^2}}\alpha_n(\Lambda_{T,\boldsymbol{\beta}})^{\frac{2R}{n}-\epsilon}d\boldsymbol{\beta}\\
            &\ll  \sup_{\frac{B}{T^2}\leq \|\boldsymbol{\beta}\|_{\infty}\leq \frac{2B}{T^2}}  (\alpha_n(\Lambda_{T,\boldsymbol{\beta}}))^{\epsilon+\frac{1}{2}-\frac{2R}{n}} \int_{\boldsymbol{\beta}\in [0,1)^R}\alpha_n(\Lambda_{T,\boldsymbol{\beta}})^{\frac{2R}{n}-\epsilon}d\boldsymbol{\beta}\\
            &\ll T^{\frac{n}{2}-2R-\eta_3},
        \end{aligned}
    \end{equation}
    for some $\eta_3>0.$ Thus, by substituting this into the last expression in $(\ref{key step for absolute convergence of the singular integrals})$, we obtain
    \begin{equation}\label{mean value over beta scaling}
            \int_{B\leq \|\boldsymbol{\beta}\|_{\infty}\leq 2B}|I(\boldsymbol{\beta};1)|d\boldsymbol{\beta}\ll B^{-\tau},
    \end{equation}
for some $\tau>0,$ where we have used $B=T^{\delta}$ with $\delta<\frac{1}{R+1}$. Furthermore, the inequality $\mathfrak{J}\ll 1$ immediately follows by $(\ref{mean value over beta scaling})$.
    \end{proof}

\bigskip

\begin{proof}[Proof of Lemma 3.2]
After decreasing $\tau>0$ if necessary, we deduce from Lemma $\ref{mean value estimate over major arcs}$, $\ref{absolute convergence of the singular series}$ and  $\ref{absolute convergence of singular integral}$ that whenever $n\geq 4R+1$, one has
\begin{equation*}
    \begin{aligned}
\int_{\mathfrak{M}_{T,\delta}}S(\bxi;T)d\bxi&=\mathfrak{S}(T^{\delta})\mathfrak{J}(T^{\delta};1)T^{n-2R}+O(T^{n-2R-\tau})\\
&=(\mathfrak{S}+(\mathfrak{S}(T^{\delta})-\mathfrak{S}))(\mathfrak{J}+ (\mathfrak{J}(T^{\delta};1)-\mathfrak{J}))T^{n-2R}+O(T^{n-2R-\tau})\\
&=\mathfrak{S}\mathfrak{J}T^{n-2R}+O(T^{n-2R-\tau}).
    \end{aligned}
\end{equation*}
It follows by \cite[Section 7]{B62} and \cite[Lemma 2.6(ii)]{{RM18}} that 
\begin{equation}\label{product of padic densities}
    \mathfrak{S}=\prod_{p:\text{prime}}\lim_{k\rightarrow \infty}\frac{1}{p^{k(n-R)}}\#\{\boldsymbol{b}\in (\mathbb{Z}/p^k\mathbb{Z})^n:\ Q_i(\boldsymbol{b})\equiv0\ (\text{mod}\ p^k)\ \text{for all }1\leq i\leq R\}.
\end{equation}
and
\begin{equation}\label{real density}
    \mathfrak{J}=\lim_{T\rightarrow \infty}
    \frac{1}{T^{n-2R}}\text{Leb}\bigg\{\mathbf{y}\in [-T,T]^n: |Q_i(\mathbf{y})|\leq \frac12 \text{for all }1\leq i\leq R\bigg\},
\end{equation}
where $\operatorname{Leb}$ denotes the Lebesgue measure on $\R^n$. 

Therefore, we conclude by \cite[Proof of Theorem 1.3]{RM18} with $d=2$ that the constant $\mathfrak{S}\mathfrak{J}$ is positive, provided that the system \eqref{systemofquadraticforms} has a non-trivial solution over $\R$ and $\Z_p$ for all primes $p$. 
\end{proof}

\bigskip

Finally, we provide the asymptotic formula for the second mean value in $(\ref{mean value of the exponential sum over major arcs})$. Before introducing the statement of this lemma, we define
\begin{equation}\label{def of I_w}
\begin{aligned}I_{\omega}(\boldsymbol{\beta};T):=\int_{\R^n}\omega\left(\frac{\boldsymbol{\gamma}}{T}\right)e\left(\sum_{i=1}^R\beta_iQ_i(\boldsymbol{\gamma})\right)d\boldsymbol{\gamma},
    \end{aligned}
\end{equation}
and
\begin{equation}\label{def of J_w}
\mathfrak{J}_{\omega}:=\int_{\boldsymbol{\beta}\in \mathbb{R}^R}I_{\omega}(\boldsymbol{\beta};1)d\boldsymbol{\beta}.
\end{equation}
\begin{lem}\label{asymptotic formula of mean value over major arcs 2}
     Let $n\geq 4R+1$ and fix $0<\delta<\frac{1}{2R+3}.$ Then one has
    \begin{equation*}
    \int_{\mathfrak{M}_{T,\delta}}S_{\omega}(\bxi;T)d\bxi=\mathfrak{S}\mathfrak{J}_{\omega}T^{n-2R}+O(T^{n-2R-\tau}),
    \end{equation*}
   for some $\tau>0,$ where the singular series $\mathfrak{S}$ is defined in \eqref{eq:singularseriesdef} depending on $Q_1,\ldots,Q_R$ and the singular integral $\mathfrak{J}_{\omega}$ is defined in \eqref{def of J_w} depending on  $Q_1,\ldots,Q_R$. Furthermore, the constant $\mathfrak{S}\mathfrak{J}_{\omega}$ is positive, provided that the system \eqref{systemofquadraticforms} has a non-trivial solution over $\R$ and $\Z_p$ for all primes $p$.
\end{lem}
\begin{proof}
We argue as in the proof of Lemma \ref{mean value estimate over major arcs}.  With the same notation in the proof of Lemma \ref{mean value estimate over major arcs}, it follows that one has
\begin{equation}\label{first step in the proof of Lemma 3.6}
\begin{aligned}
   S_{\omega}(\bxi;T)&=\sum_{\mathbf{r}\in\{0,\ldots,q-1\}^n}\sum_{\substack{\mathbf{y}\in\mathbb{Z}^n}}\omega\left(\frac{q\mathbf{y}+\mathbf{r}}{T}\right)e\bigg(\sum_{i=1}^R\bigg(\beta_i+\frac{p_i}{q}\bigg)Q_i(q\mathbf{y}+\mathbf{r})\bigg)\\
&=\sum_{\mathbf{r}\in\{0,\ldots,q-1\}^n}e\bigg(\sum_{i=1}^R\frac{p_i}{q}Q_i(\mathbf{r})\bigg)\cdot \sum_{\substack{\mathbf{y}\in\mathbb{Z}^n}}\omega\left(\frac{q\mathbf{y}+\mathbf{r}}{T}\right)e\bigg(\sum_{i=1}^R\beta_iQ_i(q\mathbf{y}+\mathbf{r})\bigg).
\end{aligned}
\end{equation}

By the standard argument using the mean value theorem, as in \cite[Lemma 5.1]{B62}, one has
    \begin{equation}\label{the standard mean value theorem2}
    \begin{aligned}
&\sum_{\substack{\mathbf{y}\in\mathbb{Z}^n}}\omega\left(\frac{q\mathbf{y}+\mathbf{r}}{T}\right)e\bigg(\sum_{i=1}^R\beta_iQ_i(q\mathbf{y}+\mathbf{r})\bigg)\\
&=q^{-n}I_{\omega}(\boldsymbol{\beta};T)+O\bigg(q^{-n}T^{n-1}\big(q+\sum_{i=1}^RqT^2|\beta_i|\big)\bigg),
    \end{aligned}
\end{equation}
where $I_{\omega}(\boldsymbol{\beta};T)$ is defined in (\ref{def of I_w}). 

Substituting (\ref{the standard mean value theorem2}) into (\ref{first step in the proof of Lemma 3.6}), we see that whenever $\bxi\in \mathfrak{M}_{T,\delta}(q,\mathbf{p}),$ one has
    \begin{equation}\label{approximation in major arcs2}
        S_{\omega}(\bxi;T)=q^{-n}S_q(\mathbf{p})I_{\omega}(\boldsymbol{\beta};T)+O(qT^{n-1+\delta}).
    \end{equation}
    Following the same argument leading from $(\ref{approximation in major arcs})$ to (\ref{last step of the proof of Lemma 3.3}), one has
 \begin{equation}\label{approximation for major arcs mean value}
\int_{\mathfrak{M}_{T,\delta}}S_{\omega}(\bxi;T)d\bxi=\mathfrak{S}(T^{\delta})\mathfrak{J}_{\omega}(T^{\delta};1)T^{n-2R}+O(T^{n-2R-\tau}),
    \end{equation}
    for some $\tau>0,$ where
\begin{equation*}
\mathfrak{J}_{\omega}(T^{\delta};1):=\int_{\|\boldsymbol{\beta}\|_{\infty}\leq T^{\delta}}I_{\omega}(\boldsymbol{\beta};1)d\boldsymbol{\beta}.
\end{equation*}

By the same argument as in the proof of Lemma \ref{absolute convergence of singular integral} with (\ref{approximation in major arcs}) and (\ref{exponential sum estimate with smooth weight}) replaced by (\ref{approximation in major arcs2}) and (\ref{exponential sum estimate with smooth weight2}), we infer that  whenever $n\geq 4R+1$, the truncated singular integral $\mathfrak{J}_{\omega}(B;1)$ converges absolutely, as $B\rightarrow \infty$. Furthermore, one has
 \begin{equation}\label{estimates for Jomega}
     |\mathfrak{J}_{\omega}-\mathfrak{J}_{\omega}(B;1)|\ll B^{-\tau}\ \ \text{and}\ \ \mathfrak{J}_{\omega}\ll 1,
 \end{equation}
 for some $\tau>0.$

After decreasing $\tau>0$ if necessary, we deduce from  $(\ref{approximation for major arcs mean value})$, $(\ref{estimates for Jomega})$, and Lemma  $\ref{absolute convergence of the singular series}$ that whenever $n\geq 4R+1$, one has
\begin{equation*}
    \begin{aligned}
\int_{\mathfrak{M}_{T,\delta}}S_{\omega}(\bxi;T)d\bxi&=\mathfrak{S}(T^{\delta})\mathfrak{J}_{\omega}(T^{\delta};1)T^{n-2R}+O(T^{n-2R-\tau})\\
&=(\mathfrak{S}+(\mathfrak{S}(T^{\delta})-\mathfrak{S}))(\mathfrak{J}_{\omega}+ (\mathfrak{J}_{\omega}(T^{\delta};1)-\mathfrak{J}_{\omega}))T^{n-2R}+O(T^{n-2R-\tau})\\
&=\mathfrak{S}\mathfrak{J}_{\omega}T^{n-2R}+O(T^{n-2R-\tau}).
    \end{aligned}
\end{equation*}
Suppose now that the system \eqref{systemofquadraticforms} has a non-trivial solution over $\mathbb{R}$ and $\mathbb{Z}_p$ for all primes $p.$ Since $\omega\equiv1$ in a neighborhood of the origin, any nonzero real solution can be rescaled to a solution at which $\omega$ is positive. This solution is nonsingular by the smoothness hypothesis. Thus the local function $\mathfrak{J}_w>0$, by \cite[Proof of Theorem 1.3]{RM18} with $d=2$. As in Lemma \ref{asymptotic formula for the mean value over major arcs}, the $p$-adic hypothesis gives $\mathfrak{S}>0$. Hence the constant $\mathfrak{S}\mathfrak{J}_{\omega}$ is positive.
\end{proof}

\bigskip

\subsection{Minor arcs analysis}
In this subsection, we provide upper bounds for the following integrals over the minor arcs:
\begin{equation*}
    \int_{\mathfrak{m}_{T,\delta}}S(\bxi;T)d\bxi \ \ \text{and}\ \ \int_{\mathfrak{m}_{T,\delta}}S_{\omega}(\bxi;T)d\bxi.
\end{equation*}
By making use of the main moment bound $(\ref{moment bound})$ together with the second part of Lemma \ref{pointwise bounds}, we have the following.
\begin{lem}\label{minor arcs estimates}
   Whenever $n\geq 4R+1,$ one has
    \begin{equation*}
         \int_{\mathfrak{m}_{T,\delta}}S(\bxi;T)d\bxi\ll T^{n-2R-\tau},
    \end{equation*}
    for some $\tau>0.$
\end{lem}
\begin{proof}
Choose $0<\lambda<\tfrac{2R}{n}$ sufficiently close to $\tfrac{2R}{n}$, and then choose
$\varepsilon>0$ sufficiently small. By Lemma~\ref{lemma2.2}, one has
\begin{equation}\label{initial bound of minor acrs}
\begin{aligned}
\int_{\mathfrak m_{T,\delta}}|S(\bxi;T)|d\bxi
&\ll T^{n/2}(\log T)^n\int_{\mathfrak m_{T,\delta}}
\alpha_n(\Lambda_{T,\bxi})^{1/2}d\bxi\\
&\ll T^{n/2}(\log T)^n
\sup_{\bxi\in\mathfrak m_{T,\delta}}
\alpha_n(\Lambda_{T,\bxi})^{1/2-\lambda}
\int_{[0,1]^R}\alpha_n(\Lambda_{T,\bxi})^\lambda d\bxi.
\end{aligned}
\end{equation}
Using Lemma~\ref{pointwise bounds} and Proposition~\ref{prop:mainmomentestimate},
the right-hand side is
\[
\ll T^{n/2}(T^{n-\delta})^{1/2-\lambda}T^\varepsilon.
\]
Since $1/2-\tfrac{2R}{n}>0$, by taking $\lambda$ close enough to $\tfrac{2R}{n}$ and then
$\varepsilon$ small enough, the exponent of $T$ is at most $n-2R-\tau$ for some
$\tau>0$. This proves the lemma.
\end{proof}

\bigskip

\begin{lem}\label{minor arcs estimates2}
   Whenever $n\geq 4R+1,$ one has
    \begin{equation*}
\int_{\mathfrak{m}_{T,\delta}}S_{\omega}(\bxi;T)d\bxi\ll T^{n-2R-\tau},
    \end{equation*}
    for some $\tau>0.$
\end{lem}
\begin{proof}
    By replacing Lemma \ref{lemma2.2} with Lemma \ref{lemma2.1} in the proof of Lemma \ref{minor arcs estimates}, we complete the proof of Lemma \ref{minor arcs estimates2}.
\end{proof}

\subsection{Proofs of the main theorems}
\begin{proof}[Proof of Theorem~\ref{thm:Meyerforsystem}, \ref{main theorem without smooth weight}, \ref{main theorem with smooth weight}] Fix major arcs parameter $0<\delta_0<1/(2R+3)$. Applying  Lemma \ref{asymptotic formula for the mean value over major arcs} and \ref{minor arcs estimates} with this parameter, and taking $\delta>0$ to be the smaller of their saving exponents,  whenever $n\geq 4R+1,$ one obtains    \begin{equation*}
\int_{[0,1]^R}S(\bxi;T)d\bxi=\mathfrak{S}\mathfrak{J}T^{n-2R}+O(T^{n-2R-\delta}),
    \end{equation*}
   for some $\delta>0,$ where the constant $\mathfrak{S}\mathfrak{J}$ is positive, provided that the system \eqref{systemofquadraticforms} has a non-trivial solution over $\R$ and $\Z_p$ for all primes $p$. This completes the proof of Theorem \ref{main theorem without smooth weight}.

  Next,  following the proof right above with Lemma \ref{asymptotic formula for the mean value over major arcs} and \ref{minor arcs estimates} replaced by Lemma \ref{asymptotic formula of mean value over major arcs 2} and \ref{minor arcs estimates2}, we complete the proof of Theorem~\ref{main theorem with smooth weight}.

     Finally, by the work of Heath-Brown \cite{HB10}, the system $(\ref{systemofquadraticforms})$ is soluble over $\mathbb{Q}_p$ for $p>(2R)^R$. Hence, by Theorem \ref{main theorem without smooth weight}, we complete the proof of Theorem \ref{thm:Meyerforsystem}.
\end{proof}

\bigskip

\section{Proof of the main moment estimate}\label{section:mainmomentestimate}

In this section, we deduce Proposition~\ref{prop:mainmomentestimate} from a uniform mean-value estimate in exterior powers, stated in Proposition~\ref{prop:localcontraction} and proved in Section~\ref{sec:sublevelestimate}. The argument follows the method of \cite{EMM98}, based on Margulis inequalities, which have been extensively studied in homogeneous dynamics. Adapting their construction of weighted sums of the height functions $\alpha_j^\lambda$, we obtain a mean-value inequality with a multiplicative constant $c>1$ independent of the averaging time $t$. Iterating this inequality at a sufficiently large fixed time yields the required $T^{o(1)}$ bound.

Throughout Sections~\ref{section:mainmomentestimate} and \ref{sec:sublevelestimate}, we assume that $\mathcal{Z}_{\mathbb C}$ is a smooth complete intersection and that $n\geq 4R+1$. We equip each exterior power $\bigwedge^k\mathbb{R}^{2n}$ with the norm induced by the Euclidean norm on $\mathbb{R}^{2n}$, and suppress the representation symbol in expressions such as $a(t)u(\boldsymbol\xi)v$.

\bigskip

\subsection{Uniform mean-value inequality in exterior powers}
For $1\leq k\leq 2n-1$, let $V_k:=\bigwedge^k\R^{2n}$. Let
$\rho_k:G\to \operatorname{GL}(V_k)$ denote the representation induced from
the standard action on $\R^{2n}$. Thus, for a decomposable vector
$v=v_1\wedge\cdots\wedge v_k$ with $v_1,\ldots,v_k\in \mathbb{R}^{2n}$, we have
$$\rho_k(g)v=(gv_1)\wedge\cdots\wedge(gv_k).$$

\begin{prop}\label{prop:localcontraction}
Assume that $\mathcal{Z}_\C$ is a smooth complete intersection and that $n\ge4R+1$.
Let $1\leq k\leq n$ and let $0<\lambda<\tfrac{2R}{n}$.
There is a constant $C_0=C_0(k,\lambda)>0$ such that for every $t\geq0$ and every nonzero $v\in V_k$,
\[
\int_{[0,1]^R}\|a(t)u(\bxi)v\|^{-\lambda}\,d\bxi
\le C_0\|v\|^{-\lambda}.
\]
\end{prop}
The proof of Proposition~\ref{prop:localcontraction} will be given in Section~5. In the remainder of this section, we deduce Proposition~\ref{prop:mainmomentestimate} from Proposition~\ref{prop:localcontraction}.

\bigskip

\subsection{Mean-value estimates for symplectic lattices}

We now deduce the mean-value inequality in Proposition~\ref{prop:globalcontraction} for lattice height functions from the estimate for individual exterior vectors in Proposition~\ref{prop:localcontraction}. Following \cite[Section 5]{EMM98}, for each fixed dimension and fixed averaging time, we distinguish two cases. If one rational subspace of the given dimension has covolume sufficiently smaller than those of all the others, it remains the minimizing subspace even after the averaging time, and Proposition~\ref{prop:localcontraction} applies directly. Otherwise, if two distinct rational subspaces have covolumes comparable to the minimum, then the intersection–sum inequality \cite[Lemma 5.2]{EMM98} bounds the corresponding height in terms of height functions of other degrees. We then combine these estimates using a suitable linear combination of the functions $\alpha_j^\lambda$.

For symplectic lattices, the symmetry \eqref{eq:symplecticsymmetry} allows us
to work only with the degrees $1,\ldots,n$. Fix $0<\lambda<\tfrac{2R}{n}$. Given an averaging time $t$, we define
\begin{equation}\label{eq:baralphadef}
    \bar{\alpha}_{\lambda,t}(\Lambda)
:=\sum_{j=1}^n \epsilon_j(t)\alpha_j(\Lambda)^\lambda,
\end{equation}
where the positive coefficients $\epsilon_1(t),\ldots,\epsilon_n(t)$ are chosen by
\begin{equation}\label{eq:height-coefficients}
\epsilon_i(t)=\exp\bigl(10n(n-i)^2t\bigr).
\end{equation}

\begin{prop}\label{prop:globalcontraction}
For every $0<\lambda<\tfrac{2R}{n}$ there is a constant
$c=c(\lambda)>1$, independent of $t$, with the following property. For every
sufficiently large $t$, we have
\begin{equation}\label{eq:globalcontraction}
\int_{[0,1]^R}\bar{\alpha}_{\lambda,t}(a(t)u(\bxi)\Lambda)\,d\bxi
\le c\bar{\alpha}_{\lambda,t}(\Lambda)+(\mathsf{D}e^t)^{20n^3}
\end{equation}
for every symplectic lattice $\Lambda\in X$, where $\mathsf{D}>1$ is a constant depending only on $Q_1,\ldots,Q_R$.
\end{prop}

The weights $\epsilon_i(t)$ and the additive constant $(\mathsf{D}e^t)^{20n^3}$ in \eqref{eq:globalcontraction} depend on the averaging time $t$, but the multiplicative constant $c$ is independent of $t$. This independence of $c$ from $t$ is the key to obtaining the arbitrarily small power loss in Proposition~\ref{prop:mainmomentestimate}.

\begin{proof}
It is enough to estimate the terms $\alpha_j^\lambda$ with $1\leq j\leq n$;
if an index larger than $n$ occurs, we use the symplectic symmetry
\eqref{eq:symplecticsymmetry}. Let $C_0$ be larger than the constants in
Proposition~\ref{prop:localcontraction} for all $1\leq j\leq n$.

There exists a constant $\mathsf{D}>1$, depending only on $Q_1,\ldots,Q_R$, such that
\begin{equation}\label{eq:vector-distortion}
(\mathsf{D}e^t)^{-1}\|w\|\leq \|a(t)u(\bxi)w\|\leq \mathsf{D}e^t\|w\|    
\end{equation}
for any $t\geq 0$, $\bxi\in [0,1]^R$, and $w\in \mathbb{R}^{2n}$. Taking exterior powers, for every $\Lambda$-rational subspace $W$ of
any dimension, it follows that
\begin{equation}\label{eq:covolume-distortion}
(\mathsf{D}e^t)^{-2n}d_\Lambda(W)
\le d_{a(t)u({\bxi})\Lambda}\big(a(t)u({\bxi}) W\big)
\le (\mathsf{D}e^t)^{2n} d_\Lambda(W).
\end{equation}

Fix $j\leq n$, and choose a $j$-dimensional $\Lambda$-rational subspace
$V$ such that $d_\Lambda(V)^{-1}=\alpha_j(\Lambda)$. The supremum in the
definition of $\alpha_j$ is attained because there are only finitely many
primitive $j$-vectors of bounded norm. Let $\mathcal P_j(\Lambda)$ denote the
finite set of $j$-dimensional $\Lambda$-rational subspaces $W$ satisfying
\[
d_\Lambda(W)\le (\mathsf{D}e^t)^{4n}d_\Lambda(V).
\]

\textbf{Case 1.} $\#\mathcal P_j(\Lambda)=1$. 

In this case we have $\mathcal P_j(\Lambda)=\{V\}$. If $W$ is a $\Lambda$-rational subspace with $\dim W=j$, minimizing $d_{a(t)u(\bxi)\Lambda}(a(t)u(\bxi) W)$, then \eqref{eq:covolume-distortion} gives
\begin{equation}
    d_\Lambda(W)
\le (\mathsf{D}e^t)^{2n} d_{a(t)u(\bxi)\Lambda}(a(t)u(\bxi) W)
\le (\mathsf{D}e^t)^{2n} d_{a(t)u(\bxi)\Lambda}(a(t)u(\bxi) V)
\le (\mathsf{D}e^t)^{4n} d_\Lambda(V),
\end{equation}
hence $W=V$. Let $v_V\in\bigwedge^j\R^{2n}$ denote the primitive Pl\"{u}cker vector representing $V\cap\Lambda$, so that $d_\Lambda(V)=\|v_V\|$. Then we have
\begin{equation}\label{eq:minimizigsubspace}
    \alpha_j(a(t)u(\bxi)\Lambda)=\|a(t)u(\bxi)v_V\|^{-1}.
\end{equation}
Thus, by Proposition~\ref{prop:localcontraction}, it follows that
\begin{equation}\label{eq:case1-local-bounded}
\int_{[0,1]^R}\alpha_j(a(t)u(\bxi)\Lambda)^\lambda d\bxi=\int_{[0,1]^R}\|a(t)u(\bxi)v_V\|^{-\lambda} d\bxi\leq C_0\|v_V\|^{-\lambda}= C_0\alpha_j(\Lambda)^\lambda.
\end{equation}

\textbf{Case 2.} $\#\mathcal P_j(\Lambda)\geq2$. 

In this case, we can find distinct
$W_1,W_2\in\mathcal P_j(\Lambda)$, and write
$$\dim(W_1\cap W_2)=j-\ell,\qquad \dim(W_1+W_2)=j+\ell$$
with $1\leq \ell\leq \min\{j,2n-j\}$. According to the intersection-sum inequality \cite[Lemma 5.2]{EMM98}, we have
\begin{equation}\label{eq:EMMintersectionsuminequality}
    d_\Lambda(W_1\cap W_2)d_\Lambda(W_1+W_2)
\le d_\Lambda(W_1)d_\Lambda(W_2).
\end{equation}
Together with the definition of
$\mathcal P_j(\Lambda)$, it gives
\begin{equation}
    \begin{aligned}
        \alpha_{j-\ell}(\Lambda) \alpha_{j+\ell}(\Lambda)&\geq d_\Lambda(W_1\cap W_2)^{-1}d_\Lambda(W_1+W_2)^{-1}\\&\geq d_\Lambda(W_1)^{-1}d_\Lambda(W_2)^{-1}\geq (\mathsf{D}e^t)^{-8n} \alpha_j(\Lambda)^2,
    \end{aligned}
\end{equation}
where $\alpha_0=\alpha_{2n}=1$. Also
$\alpha_j(a(t)u(\bxi)\Lambda)^\lambda\leq (\mathsf{D}e^t)^{n}\alpha_j(\Lambda)^\lambda$. It follows that
\begin{equation}\label{eq:case2-product-new}
\alpha_j(a(t)u(\bxi)\Lambda)^\lambda
\leq (\mathsf{D}e^t)^{3n}
\alpha_{j-\ell}(\Lambda)^{\lambda/2}
\alpha_{j+\ell}(\Lambda)^{\lambda/2}.
\end{equation}
If one of the indices $j-\ell,j+\ell$ is larger than $n$, we replace it by its symplectic complementary index, using the symmetry \eqref{eq:symplecticsymmetry}. Thus every right-hand side in \eqref{eq:case2-product-new} is bounded by one of finitely many expressions $(\mathsf{D}e^t)^{3n}
\alpha_{p}(\Lambda)^{\lambda/2}
\alpha_{q}(\Lambda)^{\lambda/2}$, where $0\leq p,q\leq n$. Explicitly, the first index $p$ is $j-\ell$, and the second index $q$ is given by
$$q=\begin{cases}
j+\ell,&j+\ell\le n,\\
2n-j-\ell,&j+\ell>n.
\end{cases}$$

If $p,q\ge1$, then in
both cases in the definition of $q$, a direct calculation gives
\begin{equation}\label{eq:coefficientsconvexity}
    \frac{\epsilon_j(t)}{(\epsilon_p(t)\epsilon_q(t))^{\frac{1}{2}}}
=e^{-10n\ell^2t}\le e^{-10nt}.
\end{equation}
Applying Young's inequality together with \eqref{eq:coefficientsconvexity}, we get
\begin{equation}\label{eq:nozeropqcase}
(\mathsf{D}e^t)^{3n}\epsilon_j(t)\alpha_p^{\lambda/2}\alpha_q^{\lambda/2}\le\epsilon_p(t)\alpha_p^{\lambda}+\epsilon_q(t)\alpha_q^{\lambda}.
\end{equation}
If $p=0$ and $q\geq 1$, similarly as above, Young's inequality gives
\begin{equation}\label{eq:zeropcase}
    (\mathsf{D}e^t)^{3n}\epsilon_j(t)\alpha_q^{\lambda/2} \le\epsilon_q(t)\alpha_q^{\lambda}+\frac{(\mathsf{D}e^t)^{6n}\epsilon_j(t)^2}{\epsilon_q(t)}.
\end{equation}
The case $q=0$ and $p\geq 1$ is similar by symmetry. It is also easy to check the case $p=q=0$. Thus, combining \eqref{eq:nozeropqcase} and \eqref{eq:zeropcase}, every product term in \eqref{eq:case2-product-new} satisfies
\begin{equation}\label{eq:case2-absorbed-new}
\epsilon_j(t)\alpha_j(a(t)u(\bxi)\Lambda)^\lambda
\le
2\bar{\alpha}_{\lambda,t}(\Lambda)+(\mathsf{D}e^t)^{10n^3}.
\end{equation}

Combining \eqref{eq:case1-local-bounded} and \eqref{eq:case2-absorbed-new},
and summing over $1\leq j\leq n$, after increasing $\mathsf{D}$ if necessary, we conclude that
$$\int_{[0,1]^R}\bar{\alpha}_{\lambda,t}(a(t)u(\bxi)\Lambda)\,d\bxi\le (C_0+2n)\bar{\alpha}_{\lambda,t}(\Lambda)+(\mathsf{D}e^t)^{20n^3},$$
as desired. Note that the number $c=C_0+2n$ is independent of $t$.
\end{proof}

\subsection{Proof of Proposition~\ref{prop:mainmomentestimate}}
We now deduce Proposition~\ref{prop:mainmomentestimate} from Proposition~\ref{prop:globalcontraction} by an iteration argument, as in \cite[Section 5]{EMM98}.

The case $\lambda=0$ is trivial, so let $0<\lambda<\tfrac{2R}{n}$, and fix
$\varepsilon>0$. Let $c>1$ be the constant in Proposition~\ref{prop:globalcontraction}.
Choose $t$ sufficiently large so that
$$\frac{\log c}{t}<\varepsilon.$$
For this fixed $t$, choose the height $\bar{\alpha}_{\lambda,t}$ as in Proposition~\ref{prop:globalcontraction}. For a nonnegative Borel function $f:X\to[0,\infty]$, define the averaging operator $\mathcal{A}_t$ by
\begin{equation}\label{eq:averagingoperatordefinition}
    \mathcal{A}_t f(\Lambda):=\int_{[0,1]^R}f(a(t)u(\bxi)\Lambda)\,d\bxi,
\end{equation}
allowing the value $+\infty$. Then Proposition~\ref{prop:globalcontraction} gives
\begin{equation}\label{eq:Margulisinequality}
    \mathcal{A}_t\bar{\alpha}_{\lambda,t}\le c\bar{\alpha}_{\lambda,t}+(\mathsf{D}e^t)^{20n^3}.
\end{equation}
Iterating \eqref{eq:Margulisinequality}, we get
\begin{equation}\label{eq:uniformboundednessbyiteration}
    (\mathcal{A}_t)^N\bar{\alpha}_{\lambda,t}(\Z^{2n})
\le c^N\left(\bar{\alpha}_{\lambda,t}(\Z^{2n})+\frac{(\mathsf{D}e^t)^{20n^3}}{c-1}\right)
\ll_t c^N .
\end{equation}
We next write the iterated average as an average over a single closed unipotent orbit. Since the forms have integer coefficients, every $A_i$ is an integral symmetric matrix, and hence $u(\mathbf n)\in\Gamma$ for any $\mathbf n\in\Z^R$. Also,
$a(t)u(\bxi)a(-t)=u(e^{2t}\bxi)$. It follows that
\begin{equation}\label{eq:randomwalkcontraction}
    a(t)u(\bxi_N)\cdots a(t)u(\bxi_1)
=a(Nt)u\bigl(\bxi_1+e^{-2t}\bxi_2+\cdots+e^{-2(N-1)t}\bxi_N\bigr).
\end{equation}
For $N\ge1$ and fixed $\bxi_2,\ldots,\bxi_N$, the function of $\bxi_1$ on the RHS of \eqref{eq:randomwalkcontraction} is
$\Z^R$-periodic after it is applied to $\Z^{2n}$. Integrating over $\bxi_1\in[0,1]^R$, we thus have
\begin{equation}\label{eq:randomwalkformulation}
(\mathcal{A}_t)^Nf(\Z^{2n})
=\int_{[0,1]^R}f(a(Nt)u(\bxi)\Z^{2n})\,d\bxi    
\end{equation}
for every nonnegative Borel function $f$. Combining \eqref{eq:uniformboundednessbyiteration} and \eqref{eq:randomwalkformulation} with
$f=\bar{\alpha}_{\lambda,t}$, we get
\begin{equation}
    \int_{[0,1]^R}
\bar{\alpha}_{\lambda,t}(a(Nt)u(\bxi)\Z^{2n})d\bxi=(\mathcal{A}_t)^N\bar{\alpha}_{\lambda,t}(\Z^{2n})
\ll_t c^N.
\end{equation}
The term $\alpha_n^\lambda$ occurs in $\bar{\alpha}_{\lambda,t}$ with coefficient $\epsilon_n(t)=1$. Therefore, when $T=e^{Nt}$ with $N\ge1$,
\begin{equation}
    \int_{[0,1]^R}
\alpha_n(\Lambda_{T,\bxi})^\lambda d\bxi
\ll_{t,\lambda} T^{\frac{\log c}{t}}
\leq T^\varepsilon .
\end{equation}
For a general $T\ge1$, write $\log T=Nt+s$, where $N\ge0$ and $0\le s<t$.
Since $\Lambda_{T,\bxi}=a(s)\Lambda_{e^{Nt},\bxi}$ and
$\alpha_n(a(s)\Lambda)\le e^{ns}\alpha_n(\Lambda)$, the same estimate holds
after changing its implied constant. The bounded range $1\le T<e^t$ is
absorbed into that constant.
This proves \eqref{moment bound}.

\bigskip

\section{Proof of Proposition~\ref{prop:localcontraction}}\label{sec:sublevelestimate}
In this section, we prove Proposition~\ref{prop:localcontraction} following the strategy of \cite[Sections 4 and 5]{KO26}. The main step is to establish uniform sublevel estimates for the polynomial maps arising from the unipotent action on exterior powers. To this end, we first derive algebraic dimension estimates for the kernels of these polynomial maps in Subsection~\ref{subsec:algebraicdimensionbound}, using the smoothness assumption on $\mathcal Z_{\mathbb C}$. We then use these dimension bounds to establish sublevel estimates first for weight vectors in Subsection~\ref{subsec:sublevelforweightvectors} and subsequently for general vectors in Subsection~\ref{subsec:sublevelforgeneralvectors}. Finally, in Subsection~\ref{subsec:proofofuniformmeanvalueestimate}, we apply these estimates to obtain the uniform mean-value bound in Proposition~\ref{prop:localcontraction}.

\bigskip

\subsection{Setting}
Let $e_1,\ldots,e_{2n}$ be the standard basis of $\R^{2n}$, and set
$$W^+=\langle e_1,\ldots,e_n\rangle,
\qquad
W^-=\langle e_{n+1},\ldots,e_{2n}\rangle.$$
Then $V_1=\mathbb R^{2n}=W^+\oplus W^-$, $\dim W^+=\dim W^-=n$, and $a(t)$ acts on $W^+$ and $W^-$ by multiplication by $e^t$ and $e^{-t}$, respectively. For $0\le \ell\le k\leq n$, set
$$V_k=\bigwedge^kV_1=\bigoplus_{\ell=0}^k V_{k,\ell},\quad V_{k,\ell}:=\left(\bigwedge^\ell W^+\right)\wedge
             \left(\bigwedge^{k-\ell}W^-\right),$$
and let \(\pi_{k,\ell}:V_k\to V_{k,\ell}\) denote the corresponding projection.

The diagonal flow
$$a(t)=
\begin{pmatrix}
e^tI_n&0\\ 0&e^{-t}I_n
\end{pmatrix}$$
acts by multiplication by $e^t$ and $e^{-t}$ on $W^+$ and $W^-$, respectively. We thus have
$$a(t)w=e^{(2\ell-k)t}w$$
for any $w\in V_{k,\ell}$. 

For $\bxi\in\R^R$, recall that $A_{\bxi}:=\sum_{i=1}^R\xi_iA_i$. With respect to the decomposition $V_1=W^+\oplus W^-$, define $N_{\bxi}:V_1\to V_1$ by
$$
N_{\bxi}=
\begin{pmatrix}
0&A_{\bxi}\\
0&0
\end{pmatrix}.
$$
Thus $N_{\bxi}$ vanishes on $W^+$ and maps $W^-$ into $W^+$ via $A_{\bxi}$, using the standard coordinate identifications $W^\pm\cong\R^n$. We extend $N_{\bxi}$ to $\bigwedge^k V_1$ as a derivation of the exterior algebra. Hence $N_{\bxi}:V_{k,\ell}\to V_{k,\ell+1}$ and $N_{\bxi}^r:V_{k,\ell}\to V_{k,\ell+r}$. Note that $N_{\mathbf{\bxi}}:V_1\to V_1$ is the infinitesimal linear map associated with the unipotent flow $s\mapsto u(s\bxi)$, hence $u(\bxi)$ acts on $\bigwedge^k V_1$ as $\exp(N_{\bxi})$. Thus, for $0\le m\le k$,
\begin{equation}\label{eq:triangular-coordinate-current}
\pi_{k,m}(u(\bxi)v)=\sum_{r=0}^{m}\frac1{r!}N_{\bxi}^r\pi_{k,m-r}(v).
\end{equation}

We shall also use the complexified spaces
$$
(V_k)_\C:=V_k\otimes_\R\C,
\qquad
(V_{k,\ell})_\C:=V_{k,\ell}\otimes_\R\C.
$$
For $\bxi\in\C^R$, the same formulas define $A_{\bxi}$ and $N_{\bxi}$ on the corresponding complexified spaces. We extend the projections $\pi_{k,\ell}$ $\C$-linearly and let $u(\bxi)$ act as $\exp(N_{\bxi})$. In particular, \eqref{eq:triangular-coordinate-current} remains valid for $\bxi\in\C^R$ and $v\in(V_k)_\C$.

\bigskip

\subsection{Dimension bounds}\label{subsec:algebraicdimensionbound}

In this subsection, we prove the algebraic dimension estimates used in the sublevel argument. For a quadratic form $Q$ on $\C^n$, let $A_Q$ denote the symmetric matrix determined by
$$
Q(\bx)=\frac12 \bx^{\intercal}A_Q\bx.
$$
Thus $\nabla Q(\bx)=A_Q\bx$. In particular, for $Q=Q_{\bxi}$ with $\bxi\in\C^R$, we write $A_Q=A_{\bxi}$ for simplicity. The following lemma provides the dimension estimates for the set $\{\bxi\in\C^R:\operatorname{corank}(A_{\bxi})\ge s\}$ for every $1\leq s\leq R.$

\begin{lem}\label{dimension estimates 1}
Assume that $\mathcal{Z}_\C$ is smooth. Then, for every $s\ge1$, one has
\begin{equation}\label{target dimension estimate}
    \dim\{\bxi\in\C^R:\operatorname{corank}(A_{\bxi})\ge s\}
\le \max\{R-s,0\}.
\end{equation}
Equivalently, for $1\le s< R$, the corresponding projective locus in $\bP^{R-1}(\C)$ has dimension at most $R-s-1$, and for $s\geq R$ it is empty.
\end{lem}
\begin{proof}
Let $\mathcal Q=\operatorname{span}_\C(Q_1,\ldots,Q_R)$. Since $\mathcal{Z}_\C$ is a complete intersection of codimension $R$, the quadratic forms $Q_1,\ldots,Q_R$ are linearly independent. Thus, $\dim\mathcal Q=R$, and we may identify $\mathcal Q$ with $\C^R$ through $\bxi\mapsto A_{\bxi}$. Hence, if we define the set $D_s:=\{[Q]\in \mathbb{P}(\mathcal{Q}):\operatorname{corank}(A_Q)\ge s\},$ we see that \begin{equation}\label{dual dimension}
     \dim\{\bxi\in\C^R:\operatorname{corank}(A_{\bxi})\ge s\}=\text{dim} (D_s)+1.
\end{equation}

 By \cite[Lemma 3.1]{BDHB15}, there is an optimal system of quadratic forms $\widetilde{Q}_1=\cdots=\widetilde{Q}_R=0$ with $\mathcal{Z}_\C=V(\widetilde{Q}_1,\ldots,\widetilde{Q}_R)$ such that for every $1\leq s\leq R$, the variety $$V(\widetilde{Q}_1,\ldots,\widetilde{Q}_s)\subseteq \mathbb{P}^{n-1}$$
is a smooth complete intersection of codimension $s$, and each of $\widetilde{Q}_i$ is a linear combination of $Q_1,\ldots, Q_R$ over $\mathbb{C}.$ Define the set $\mathcal Q_s:=\operatorname{span}_\C(\widetilde{Q}_1,\ldots,\widetilde{Q}_s)$. We notice here that $\text{dim}(\mathbb{P}(\mathcal{Q}_s))=s-1$.

We may assume that $D_s\cap \mathbb{P}(\mathcal{Q}_s)\neq\emptyset.$ Indeed, otherwise, we find by the projective dimension theorem that
\begin{equation*}
    \text{dim}(D_s) + \text{dim}(\mathbb{P}(\mathcal{Q}_s))< R-1,
\end{equation*}
which implies that $ \text{dim} (D_s)\leq R-s-1.$ This proves $(\ref{target dimension estimate})$, by $(\ref{dual dimension}).$ Let $Q$ be a quadratic form  with $[Q]\in D_s\cap \mathbb{P}(\mathcal{Q}_s)$ and consider $K:= \text{ker }A_Q$. Meanwhile, since $Q\in \mathcal{Q}_s$, there exist quadratic forms $\widetilde{Q}'_1,\ldots,\widetilde{Q}'_{s-1}$ in $\mathcal{Q}_s$ such that 
$$V(\widetilde{Q}_1,\ldots,\widetilde{Q}_s)=V(Q,\widetilde{Q}'_1,\ldots,\widetilde{Q}'_{s-1}).$$ 
Thus, since this quadratic form $Q$ vanishes at $K$, we see that 
$$V(\widetilde{Q}_1,\ldots,\widetilde{Q}_s) \cap \mathbb{P}(K)=V(\widetilde{Q}'_1,\ldots,\widetilde{Q}'_{s-1}) \cap \mathbb{P}(K).$$ Therefore,
since $\text{dim}(\mathbb{P}(K))\geq s-1$, the variety $V(\widetilde{Q}_1,\ldots,\widetilde{Q}_s) \cap \mathbb{P}(K)$ is non-empty. Meanwhile, since
$V(\widetilde{Q}_1,\ldots,\widetilde{Q}_s)\subseteq \mathbb{P}^{n-1}$ is smooth, it follows that for any quadratic form $\widetilde{Q}$ with $[\widetilde{Q}]\in \mathbb{P}(\mathcal{Q}_s)$, one has $\nabla \widetilde{Q}(\bx)\neq \boldsymbol{0}$ for all $\bx\in V(\widetilde{Q}_1,\ldots,\widetilde{Q}_s)$. However, we have $\nabla Q(\bx)=\boldsymbol{0}$ on the non-empty set $V(\widetilde{Q}_1,\ldots,\widetilde{Q}_s) \cap \mathbb{P}(K)$. This is a contradiction. Thus we have $D_s\cap \mathbb{P}(\mathcal{Q}_s)=\emptyset$, and the preceding dimension argument proves the assertion for $1\leq s\leq R.$ In particular, one has $D_R=\emptyset$. For $s>R,$ we have $D_s\subseteq D_R$, so the corresponding affine locus is contained in $\{0\}$ and has dimension at most zero. This proves the remaining cases.
\end{proof}

\begin{lem}[Injectivity of $N_{\bxi}^r$]
\label{lem:injectivity-Nr}
Let $1\leq k \leq n$, $0\leq \ell\leq k$, and \(0\le r\le k-\ell\). For any $\bxi\in\mathbb{R}^R$ with
\(\operatorname{rank}(A_{\bxi})\ge n-k+2\ell+r\), the map \(N_{\bxi}^r:V_{k,\ell}\to V_{k,\ell+r}\) is injective. The same assertion holds for $\bxi\in\C^R$, with $V_{k,\ell}$ and $V_{k,\ell+r}$ replaced by their complexifications.
\end{lem}

\begin{proof}
The case $r=0$ is immediate. Set
$s=\operatorname{rank}(A_{\bxi})$ and $S=\{1,\ldots,s\}$. Choose bases
$e_1^+,\ldots,e_n^+$ of $W^+$ and
$e_1^-,\ldots,e_n^-$ of $W^-$ such that
$$
N_{\bxi}e_i^-=
\begin{cases}
e_i^+,&1\le i\le s,\\
0,&i>s.
\end{cases}
$$
For $I,J\subseteq \{1,\ldots,n\}$ we will denote $$e_{I}^{+}:=\bigwedge_{i\in I}e_i^+, \quad e_J^{-}:=\bigwedge_{j\in J}e_j^-.$$ The vectors $e_I^+\wedge e_J^-$, where $I,J\subseteq\{1,\ldots,n\}$ and $|I|+|J|=k$, form a basis of $V_k$. For each such pair $(I,J)$, define
$$
\omega_{I,J}=
\left(\bigwedge_{i\in I\setminus S}e_i^+\right)
\wedge
\left(\bigwedge_{j\in J\setminus S}e_j^-\right)
\wedge
\left(\bigwedge_{i\in I\cap J\cap S}
      (e_i^+\wedge e_i^-)\right).
$$
Observe that $N_{\bxi}\omega_{I,J}=0$. Let $L=\bigl((I\cup J)\setminus(I\cap J)\bigr)\cap S$. For $F\subseteq L$ and $i\in L$, let
$$
v_F=\omega_{I,J}\wedge\bigwedge_{i\in L}w_{i,F},
\qquad
w_{i,F}=
\begin{cases}
e_i^+,&i\in F,\\
e_i^-,&i\in L\setminus F,
\end{cases}
$$
where the factors within each indexed wedge product are
ordered by increasing index. Define
$$
E_{I,J}:=\operatorname{span}\{v_F:F\subseteq L\}
\subseteq V_k.
$$
Each $v_F$ is, up to sign, a standard basis vector of $V_k$.
The distinct subspaces $E_{I,J}$ are spanned by disjoint
subsets of this basis and give a direct sum decomposition
of $V_k$. Moreover, since $N_{\bxi}\omega_{I,J}=0$, we have
$$
N_{\bxi}v_F
=\sum_{i\in L\setminus F}v_{F\cup\{i\}}.
$$
Thus each $E_{I,J}$ is $N_{\bxi}$-invariant. It suffices
to prove injectivity on $E_{I,J}\cap V_{k,\ell}$ for
each summand whose intersection with $V_{k,\ell}$ is nonzero.

Fix such a summand $E$. We may choose $I,J$ such that
$E=E_{I,J}$, $|I|=\ell$ and $|J|=k-\ell$. Let us write $ p=|I\cap L|$ and $q=|J\cap L|$ for simplicity. Then $p+q=|L|$, and
$$
\begin{aligned}
q-p
&=|J\cap S|-|I\cap S|\\
&=k-2\ell-\bigl(|J\setminus S|-|I\setminus S|\bigr)\ge k-2\ell-(n-s)
\ge r.
\end{aligned}
$$
In particular, $2p+r\le|L|$. The vectors $v_F$ with
$|F|=p$ form a basis of $E\cap V_{k,\ell}$, and
\begin{equation}\label{eq:Nr-inclusion-map}
N_{\bxi}^r v_F
=r!\sum_{\substack{F\subseteq F'\subseteq L\\|F'|=p+r}}v_{F'}.
\end{equation}

To prove the injectivity of $N_{\bxi}^r$, suppose $N_{\bxi}^r v=0$ and write $$v=\sum_{\substack{F\subseteq L\\|F|=p}}a_Fv_F,$$
where $a_F\in\mathbb{R}$. Then in view of \eqref{eq:Nr-inclusion-map}, for every $(p+r)$-element subset $F'\subseteq L$, we have
\begin{equation}\label{eq:subset-coefficient-relation}
\sum_{\substack{F\subseteq F'\\|F|=p}}a_F=0.
\end{equation}

We claim that the relations \eqref{eq:subset-coefficient-relation} imply $a_F=0$ for every $F$
whenever $2p+r\le|L|$. We prove this claim by induction on $p$.
The case $p=0$ is immediate. Suppose that $p\ge1$.
For distinct $x,y\in L$ and
$L'\subseteq L\setminus\{x,y\}$ with $|L'|=p+r-1$,
subtracting \eqref{eq:subset-coefficient-relation} for
$L'\cup\{x\}$ and $L'\cup\{y\}$ gives
$$\sum_{\substack{D\subseteq L'\\|D|=p-1}}
\bigl(a_{D\cup\{x\}}-a_{D\cup\{y\}}\bigr)=0.$$
Since $|L\setminus\{x,y\}|\ge2(p-1)+r$, the induction hypothesis implies $a_{D\cup\{x\}}=a_{D\cup\{y\}}$ for every $(p-1)$-element subset $D\subseteq L\setminus\{x,y\}$.
Any two $p$-element subsets of $L$ can be connected by
successively replacing one element, so all coefficients
$a_F$ have a common value $a$. It follows from \eqref{eq:subset-coefficient-relation} that $\binom{p+r}{p}a=0$, hence $a=0$. This proves the claim.

We get $v=0$ immediately from the claim. Thus $N_{\bxi}^r$ is injective on each $E_{I,J}\cap V_{k,\ell}$, and therefore on $V_{k,\ell}$.
The same proof applies over $\C$.
\end{proof}

We now record the dimension estimate required for the proof of Proposition~\ref{prop:localcontraction}. For $1\leq k\leq n$, $0\leq \ell \leq k$, $1\leq r\leq k-\ell$, and non-zero $w\in (V_{k,\ell})_{\mathbb{C}}$, we define
$$Z_r(w):=\{\bxi\in\C^R:N_{\bxi}^r w=0\}.$$

\begin{lem}[Dimension estimate for $Z_r(w)$]\label{lem:zerolocusdimensionestimate}
Assume that $\mathcal{Z}_\C$ is smooth. Let $1\leq k\leq n$, $0\leq \ell<\frac{k}{2}$, and $1\leq r\leq \lceil\tfrac{k}{2}\rceil-\ell$. Then for any nonzero $w\in  (V_{k,\ell})_\C$, we have
\begin{equation}\label{eq:codimensionhigherpower}
    \operatorname{codim}_{\C^R} Z_r(w)\ge\min\{R,\,k-2\ell-r+1\}.
\end{equation}
In particular, the polynomial map $\bxi\mapsto N_{\bxi}^r w$ is not identically zero.
\end{lem}

\begin{proof}
If \(\bxi\in Z_r(w)\), then the map $N_{\bxi}^r: (V_{k,\ell})_\C\to (V_{k,\ell+r})_\C$ is not injective.  By Lemma~\ref{lem:injectivity-Nr}, this implies $\operatorname{rank}(A_{\bxi})<n-k+2\ell+r$, i.e. 
$$\operatorname{corank}(A_{\bxi})\ge k-2\ell-r+1.$$
Thus $Z_r(w)$ is contained in the locus of corank at least $k-2\ell-r+1$, hence Lemma~\ref{dimension estimates 1} gives the bound \eqref{eq:codimensionhigherpower}. Moreover, since $k-2\ell-r+1\geq 1$, we have $\operatorname{codim}_{\C^R} Z_r(w)\ge 1$, hence $\bxi\mapsto N_{\bxi}^r w$ is not identically zero. 
\end{proof}

\bigskip

\subsection{Sublevel estimates for weight vectors}\label{subsec:sublevelforweightvectors}
Let $E$ be a finite-dimensional real Euclidean space, and equip
$E_\C:=E\otimes_\R\C$ with the induced Hermitian norm. For a polynomial map
$P$ from a finite-dimensional complex Hermitian space to $E_\C$, define
its complex singularity exponent at a point $z$ of its domain by
\begin{equation}\label{eq:complexsingularityexpdef}
    c_z(P):=\sup\{s>0:\|P\|^{-2s}\text{ is locally integrable near }z\}.
\end{equation}
Integrability is with respect to Lebesgue measure on the underlying real
vector space.
We set $c_z(P)=0$ if $P\equiv0$; if $P(z)\ne0$, then $c_z(P)=+\infty$.

The following lemma records two stability statements from \cite[Appendix~A]{KO26}, in the form needed here. These statements are based on results of \cite{DemaillyKollar2001} and \cite{Cutkosky2017} in real and complex analytic geometry.
Part~(i) compares the complex singularity exponent of a polynomial map with that of its first nonzero homogeneous Taylor term, while part~(ii) gives uniform real negative-moment bounds from a common lower bound for the complex singularity exponents.

\begin{lem}\label{lem:polynomial-integrability}
\begin{enumerate}[label=\textnormal{(\roman*)}]
\item Let $P:\C^R\to E_\C$ be a nonzero polynomial map, and let $z$ be a
zero of $P$. If $H_{P,z}$ is the first nonzero homogeneous term in the
Taylor expansion of $P(z+\bxi)$, then
\begin{equation}\label{eq:homogeneousTaylorsemicontinuity}
    c_z(P)\ge c_0(H_{P,z}).
\end{equation}

\item Let $d\ge1$, let $Y\subset\R^d$ be compact, and let
$\mathcal P:\R^d\times\R^R\to E$ be polynomial. Write
$P_y(\bxi):=\mathcal P(y,\bxi)$, and let $P_y^\C$ denote its
complexification. Let $\Omega\subset\R^R$ be bounded and open. Suppose that
there are $\alpha>0$ and an open set $\mathcal{O}\subset\C^R$ containing
$\overline \Omega$ such that $c_z(P_y^\C)\ge\alpha$ for all $y\in Y$ and $z\in \mathcal{O}$. Then, for every $0<\beta<\alpha$,
\begin{equation}\label{eq:lowersemicontinuityofsingularity}
    \sup_{y\in Y}\int_\Omega\|P_y(\bxi)\|^{-\beta}\,d\bxi<\infty.
\end{equation}
\end{enumerate}
\end{lem}

\begin{proof}
For part~(i), let $h=\deg H_{P,z}$. The map
$(t,\bxi)\mapsto t^{-h}P(z+t\bxi)$ extends polynomially to $t=0$ with
value $H_{P,z}(\bxi)$. For $t\ne0$, its complex singularity exponent at
$\bxi=0$ equals $c_z(P)$. Thus, the lower semicontinuity in
\cite[Lemma~A.2(1)]{KO26} gives part~(i).
Part~(ii) follows from \cite[Proposition~A.1]{KO26}, applied to the polynomial family $(y,\bxi)\mapsto \mathcal{P}(y_0+y,\bxi)$ for each $y_0\in Y$, together with a finite covering of the compact set $Y$.
\end{proof}

For the remainder of this section, fix an integer
$1\le k\le n$ and set
$$m:=\left\lceil\frac{k}{2}\right\rceil,
\qquad
\gamma_{k,\ell}:=\min\left\{1,\frac{R}{m-\ell}\right\}$$
for $0\leq \ell<m$. We also write $\gamma_k:=\gamma_{k,0}$ and let
$\operatorname{Leb}$ denote Lebesgue measure on $\R^R$.

\begin{lem}[Sublevel estimate for weight vectors]
\label{lem:weight-vector-sublevel}
Assume that $\mathcal{Z}_\C$ is smooth and that $n\ge2R$.
Let $0\leq \ell<m$ and $0\ne w\in(V_{k,\ell})_\C$. Then we have
\begin{equation}\label{eq:complexsingularityexponent}
    c_z\bigl(\bxi\mapsto N_{\bxi}^{m-\ell} w\bigr)\ge\gamma_{k,\ell}
\end{equation}
for any $z\in\C^R$.
For $0\neq w\in V_{k,\ell}$, every bounded open \(\Omega\subset\R^R\) and every \(0<\beta<\gamma_{k,\ell}\), we also have
\begin{equation}\label{eq:realsingularityexponent}
    \int_\Omega \|N_{\bxi}^{m-\ell} w\|^{-\beta}\,d\bxi<\infty.
\end{equation}
Consequently, for every \(\eta>0\) and every
\(0<\varepsilon<1\),
\begin{equation}\label{eq:sublevelestimateoverreal}
    \operatorname{Leb}\{\bxi\in\Omega:\|N_{\bxi}^{m-\ell} w\|\le\varepsilon\|w\|\}
\ll_{w,\Omega,\eta}
\varepsilon^{\gamma_{k,\ell}-\eta}.
\end{equation}
\end{lem}

\begin{proof}
We prove the complex singularity-exponent bound \eqref{eq:complexsingularityexponent} by induction on
$r=m-\ell\ge1$. Let $F:\mathbb{C}^R\to (V_{k,m})_{\C}$ be the polynomial map defined by $F(\bxi):=N_{\bxi}^r w$, and recall the definition $Z_1(w)=\{\bxi\in\C^R:N_{\bxi} w=0\}$. Since $\bxi\mapsto N_{\bxi}w$ is complex linear,
$Z_1(w)$ is a complex linear subspace. By Lemma~\ref{lem:zerolocusdimensionestimate}, we also have $F\not\equiv0$.

It suffices to consider the zeros of $F$. Let $z$ be a zero of $F$.

\textbf{Case 1}. $z\notin Z_1(w)$. \\
This case is empty if $r=1$. Otherwise, let $j$ be the largest integer with $1\le j\le r-1$ such that
$N_z^jw\neq 0$. We write $w'=N_z^jw$ for simplicity. Since the operators $N_{\bxi}$ commute,
$$F(z+\bxi)=\sum_{i=0}^r\binom{r}{i} N_{\bxi}^iN_z^{r-i}w.$$
The terms of degree less than $r-j$ vanish. The degree-$(r-j)$ term is
$\binom {r}{r-j} N_{\bxi}^{r-j}w'$, which is not identically zero by Lemma~\ref{lem:zerolocusdimensionestimate}. Since $0\ne w'\in(V_{k,\ell+j})_\C$ and $r-j<r$, the induction hypothesis
together with Lemma~\ref{lem:polynomial-integrability}(i) gives
$$c_z(F)\ge c_0\bigl(\bxi\mapsto N_{\bxi}^{r-j}w'\bigr)
\ge\min\left\{1,\frac {R}{r-j}\right\}
\ge\gamma_{k,\ell}.$$

\textbf{Case 2}. $z\in Z_1(w)$. \\
Since the operators $N_{\bxi}$ commute, we have
\begin{equation}\label{eq:translationinvariance}
    F(x+y)=\sum_{i=0}^r\binom{r}{i} N_{y}^iN_x^{r-i}w=N_{y}^rw=F(y)
\end{equation}
for any $x\in Z_1(w)$ and $y\in\mathbb{C}^R$. In particular, $c_z(F)=c_0(F)$. 

Choose the Hermitian orthogonal complement $Z_1(w)^\perp\subset\C^R$ and let
$d:=\dim_\C Z_1(w)^\perp$. The translation invariance \eqref{eq:translationinvariance} along $Z_1(w)$ and Fubini's theorem give $$c_y(F)=c_y(F|_{Z_1(w)^\perp})$$ for any $y\in Z_1(w)^\perp$. Applying Lemma~\ref{lem:zerolocusdimensionestimate} with $r=1$, we get $d\ge\min\{R,k-2\ell\}$. Since $r=\lceil\frac{k-2\ell}{2}\rceil\le k-2\ell$, it follows that
\begin{equation}\label{eq:gammaupperbound}
    \frac{d}{r}\ge\min\left\{1,\frac Rr\right\}=\gamma_{k,\ell}.
\end{equation}

Fix $0<\beta<\gamma_{k,\ell}$, and let $0\neq y\in Z_1(w)^{\perp}$. Then $N_y w\neq 0$. If $F(y)=0$, then \textbf{Case 1} gives $c_y(F)\geq \gamma_{k,\ell}$. If $F(y)\neq 0$, then we immediately have $c_y(F)=+\infty$. It follows that $$c_y(F|_{Z_1(w)^{\perp}})=c_y(F)\geq \gamma_{k,\ell},$$ hence $\|F|_{Z_1(w)^\perp}\|^{-2\beta}$ is locally integrable near any non-zero point $y\in Z_1(w)^{\perp}$. Let $\mathcal{A}$ denote the compact annulus $\{y\in Z_1(w)^\perp:1\le\|y\|\le2\}$. By the compactness of this annulus, it follows that
\begin{equation}\label{eq:finitenessonannulus}
    \int_{\mathcal{A}}\|F|_{Z_1(w)^\perp}(y)\|^{-2\beta}\,dy<\infty.
\end{equation}
Let $\mathcal{S}$ be the unit sphere in $Z_1(w)^\perp$, with surface measure $d\sigma$.
Since $F$ is homogeneous of degree $r$, polar coordinates in the real $2d$-dimensional space $Z_1(w)^\perp$ and \eqref{eq:finitenessonannulus} give
$$\int_{\mathcal{S}}\|F(\theta)\|^{-2\beta}\,d\sigma(\theta)<\infty.$$
Again using polar coordinates, this in turn implies that
\begin{equation}\label{eq:finitenessonorthogonal}
    \int_{\substack{y\in Z_1(w)^\perp\\\|y\|\leq1}}
\|F(y)\|^{-2\beta}\,dy =
\int_0^1\rho^{2d-1-2r\beta}\,d\rho
\int_{\mathcal{S}}\|F(\theta)\|^{-2\beta}\,d\sigma(\theta)
<\infty.
\end{equation}
Since $\beta<\gamma_{k,\ell}\le \frac{d}{r}$ was arbitrary,
we conclude that $c_0(F)=c_0(F|_{Z_1(w)^{\perp}})\ge\gamma_{k,\ell}$. By the translation invariance \eqref{eq:translationinvariance}, the same bound holds at every point of $Z_1(w)$. This completes the induction, hence proves \eqref{eq:complexsingularityexponent}.

For real $w$, Lemma~\ref{lem:polynomial-integrability} (ii), applied to the constant family $P_y(\bxi)=N_{\bxi}^rw$ with $Y=\{0\}$, proves \eqref{eq:realsingularityexponent}. Choosing $\beta$ so that $\max\{0,\gamma_{k,\ell}-\eta\}<\beta<\gamma_{k,\ell}$, and applying Chebyshev's inequality, the sublevel estimate \eqref{eq:sublevelestimateoverreal} follows.

\end{proof}

\subsection{From weight vectors to general vectors}\label{subsec:sublevelforgeneralvectors}

In the remainder of this section, we assume that $\mathcal{Z}_\C$ is a smooth complete intersection and that $n\ge 4R+1$. 

Since \(n\ge4R+1\), we have \(\gamma_k\ge \frac{2R}{n}\), except when $k=n$
and $n$ is odd. Indeed, in the other cases $m\le \frac{n}{2}$, so both $1$ and
$\frac{R}{m}$ are at least $\frac{2R}{n}$. If $k=n$ and $n$ is odd, then
$m=\frac{n+1}{2}$ and $\gamma_k=\frac{2R}{n+1}$.

\begin{prop}\label{prop:triangular-polynomial-local-integrability}
Let $v_j\in(V_{k,j})_\C$ for $0\leq j\leq m$, let
\begin{equation}\label{eq:triangularpolynomialdef}
    F(\bxi)=\sum_{i=0}^{m}\frac1{i!}N_{\bxi}^i v_{m-i},
    \end{equation}
and assume that at least one of $v_0,\ldots,v_m$ is nonzero. Then for any $z\in\C^R$ we have
\begin{equation}\label{eq:singularityexponentboundfortriangular}c_z(F)\ge\gamma_k.
\end{equation}
If all $v_j$ are real, then $\|F\|^{-\beta}$ is locally integrable on
$\R^R$ for every $0<\beta<\gamma_k$.
\end{prop}

\begin{proof}
Let $z\in\C^R$ be a zero of $F$. Since
$N_{z+\bxi}=N_z+N_{\bxi}$ and the operators commute,
$$F(z+\bxi)=\sum_{i=0}^{m}\frac1{i!}N_{\bxi}^i w_{m-i},
\qquad
w_j=\sum_{i=0}^{j}\frac1{(j-i)!}N_z^{j-i}v_i\in(V_{k,j})_\C$$
for $0\leq j\leq m$. Since at least one of $v_0,\ldots,v_m$ is nonzero and the change of variables $(v_0,\ldots,v_m)\mapsto (w_0,\ldots,w_m)$ is invertible, at least one of $w_0,\ldots,w_m$ is nonzero as well. Since $F(z)=0$, we have $w_m=0$. Let $i_0$ be the smallest integer with $1\le i_0\le m$ and $w_{m-i_0}\ne0$. Then Lemma~\ref{lem:zerolocusdimensionestimate}, applied with $m-i_0<\frac{k}{2}$ and $i_0=m-(m-i_0)$, shows that $\bxi\mapsto N_{\bxi}^{i_0}w_{m-i_0}$ is not identically zero. Hence, the first nonzero homogeneous Taylor term of $F$ at $z$ is
$$H_{F,z}(\bxi)=\frac1{i_0!}N_{\bxi}^{i_0} w_{m-i_0}.$$
Lemmas~\ref{lem:polynomial-integrability} (i)
and~\ref{lem:weight-vector-sublevel} therefore give
$$c_z(F)\ge c_0(H_{F,z})
\ge\min\left\{1,\frac{R}{i_0}\right\}
\ge\gamma_k.$$
If $F(z)\neq 0$, then $c_z(F)=+\infty$.
For real coefficients, Lemma~\ref{lem:polynomial-integrability} (ii),
applied to the constant family $P_y=F$ with $Y=\{0\}$ on any bounded
open ball, proves the final assertion.
\end{proof}

For $v\in V_k$ we define
$$M(v):=\max_{0\le i\le m}\|\pi_{k,i}(v)\|.$$

\begin{prop}[Sublevel estimate for general vectors]
\label{prop:general-vector-sublevel}
Let \(\Omega\subset\R^R\) be bounded and open. For every \(\eta>0\),
\[
\operatorname{Leb}\left\{\bxi\in\Omega:
\|\pi_{k,m}(u(\bxi)v)\|\le\varepsilon M(v)\right\}
\ll_{k,\Omega,\eta}
\varepsilon^{\gamma_k-\eta}
\]
for every \(v\in V_k\) with \(M(v)>0\) and every \(0<\varepsilon<1\). In particular, if \(n\ge4R+1\) and either \(k<n\), or \(k=n\) with \(n\) even, the exponent may be replaced by \(\tfrac{2R}{n}-\eta\).
\end{prop}

\begin{proof}
Choose a bounded open Euclidean ball $\mathcal{O}$ with
\(\overline\Omega\subset \mathcal{O}\).
By homogeneity we may assume $M(v)=1$, and write $v_i=\pi_{k,i}(v)$. Then
\begin{equation}
    F_v(\bxi):=\pi_{k,m}(u(\bxi)v)
=\sum_{i=0}^{m}\frac1{i!}N_{\bxi}^{i} v_{m-i}.
\end{equation}
The tuples $(v_0,\ldots,v_m)\in\bigoplus_{i=0}^m V_{k,i}$ with
$\max_i\|v_i\|=1$ form a compact set, and the map
$(v_0,\ldots,v_m,\bxi)\mapsto F_v(\bxi)$ is polynomial.

No polynomial in the family $\{F_v: M(v)=1\}$ is identically zero. Indeed, if $i_0$ is the
smallest index with $v_{i_0}\ne0$, then the highest-degree homogeneous term
of $F_v$ is $\frac1{(m-i_0)!}N_{\bxi}^{m-i_0}v_{i_0}$. If $i_0=m$, this is the nonzero constant $v_m$. If $i_0<m$, it is not
identically zero by Lemma~\ref{lem:zerolocusdimensionestimate}.

By Proposition~\ref{prop:triangular-polynomial-local-integrability}, we have \(c_z(F_v^\C)\ge\gamma_k\) for any \(z\in\C^R\), where $F_v^\C$ denotes the complexification of $F_v$.
After choosing real coordinates on $\bigoplus_{i=0}^m V_{k,i}$,
Lemma~\ref{lem:polynomial-integrability}(ii), applied to the compact set of normalized tuples $(v_0,\ldots,v_m)$ above, gives
\begin{equation}
    \sup_{M(v)=1}\int_{\Omega}\|F_v(\bxi)\|^{-\beta}\,d\bxi<\infty
\end{equation}
for any $0<\beta<\gamma_k$. It follows from Chebyshev's inequality that
\begin{equation}\label{eq:sublevelestimatebyChebyshev}
\operatorname{Leb}\{\bxi\in \Omega:\|F_v(\bxi)\|\le\varepsilon\}
\ll_{\beta,\Omega}\varepsilon^\beta    
\end{equation}
uniformly for any $v$ with $M(v)=1$. Indeed, if $\eta<\gamma_k$, choose
$\gamma_k-\eta<\beta<\gamma_k$; if $\eta\ge\gamma_k$, \eqref{eq:sublevelestimatebyChebyshev} follows from $\operatorname{Leb}(\Omega)<\infty$ and $\varepsilon^{\gamma_k-\eta}\ge1$.
Rescaling by $M(v)$ proves the first assertion. The final assertion follows from $\gamma_k\ge\tfrac{2R}{n}$ in the stated range.
\end{proof}

\bigskip

\subsection{Proof of Proposition~\ref{prop:localcontraction}}\label{subsec:proofofuniformmeanvalueestimate}

Recall that $n\ge4R+1$, $1\le k\le n$, $m=\left\lceil\frac{k}{2}\right\rceil$, and $0<\lambda<\tfrac{2R}{n}$. Let $\mu_m=2m-k\in\{0,1\}$, so that $a(t)$ acts on $V_{k,m}$ by multiplication by $e^{\mu_m t}$. Fix $0\ne v\in V_k$.
By equivalence of norms, we may work on $V_k$ with the norm
$$\|v\|=\max_{0\le j\le k}\|\pi_{k,j}(v)\|,$$
where the component norms are Euclidean. By compactness of $[0,1]^R$ and
\eqref{eq:triangular-coordinate-current}, there is a constant $C\ge1$ such that
\[
C^{-1}\|v\|\le\|u(\bxi)v\|\le C\|v\|,
\qquad M(u(\bxi)v)\le CM(v)
\]
for all $v$ and $\bxi\in[0,1]^R$.

\textbf{Case 1}. $M(v)<(2C^2)^{-1}\|v\|$\\
If $M(v)<(2C^2)^{-1}\|v\|$, then $M(u(\bxi)v)<(2C)^{-1}\|v\|$, whereas $\|u(\bxi)v\|\ge C^{-1}\|v\|$. We observe that $m<k$ and
$$\max_{m<j\le k}\|\pi_{k,j}(u(\bxi)v)\|\ge C^{-1}\|v\|.$$
Note that for $m<j\leq k$, the $a(1)$-weight of $V_{k,j}$ is $2j-k\geq \mu_m+2\ge2$. We thus have
$$\int_{[0,1]^R}\|a(t)u(\bxi)v\|^{-\lambda}\,d\bxi
\ll e^{-2\lambda t}\|v\|^{-\lambda}.$$

\textbf{Case 2}. $M(v)\geq(2C^2)^{-1}\|v\|$\\
In this case we define
$$F(\bxi):=\frac{\|\pi_{k,m}(u(\bxi)v)\|}{M(v)},$$
$$\Xi(\varepsilon):=\{\bxi\in[0,1]^R: F(\bxi)\leq \varepsilon\}$$ 
for $\varepsilon>0$. For every $0<\gamma_0<\gamma_k$, applying Proposition~\ref{prop:general-vector-sublevel} with the open box $\Omega=(-1,2)^R$, we have
\begin{equation}\label{eq:sublevelestimateappication}
    \operatorname{Leb}\big(\Xi(\varepsilon)\big)
\ll_{\gamma_0}\min\{1,\varepsilon^{\gamma_0}\}
\end{equation}
for any $\varepsilon>0$, uniformly in $v$.
Since $\|a(t)w\| \geq e^{-kt}\|w\|$ for any $w\in V_k$, we have
\begin{equation}\label{eq:smallestweightcomponent}
    \|a(t)u(\bxi)v\|
\ge e^{-kt}\|u(\bxi)v\|
\ge C^{-1}M(v)e^{-kt}.
\end{equation}
On the other hand, we also find
\begin{equation}\label{eq:smallestnonnegativeweightcomponent}
    \|a(t)u(\bxi)v\|\ge\|\pi_{k,m}(a(t)u(\bxi)v)\|=  M(v)e^{\mu_m t}F(\bxi).
\end{equation}
Combining \eqref{eq:smallestweightcomponent} and \eqref{eq:smallestnonnegativeweightcomponent}, it follows that
\begin{equation}\label{eq:local-two-component-lower-bound}
\|a(t)u(\bxi)v\|
\ge C^{-1}M(v)e^{-kt}\max\{1,e^{(k+\mu_m)t}F(\bxi)\}.
\end{equation}

\textbf{Subcase 2-1}. $k$ is even.\\
In this subcase, $\mu_m=0$ and $m=\frac{k}{2}\le \frac{n}{2}$, so $\gamma_k\ge \frac{2R}{n}>\lambda$. Choose $\gamma_0$ so that $\lambda<\gamma_0<\gamma_k$. By
\eqref{eq:local-two-component-lower-bound}, we have
$$\|a(t)u(\bxi)v\|\ge C^{-1}M(v)\max\{e^{-kt},F(\bxi)\}.$$
In combination with \eqref{eq:sublevelestimateappication}, it follows that
\begin{equation*}
    \begin{aligned}
        \int_{[0,1]^R}\|a(t)u(\bxi)v\|^{-\lambda}\,d\bxi &\leq  \int_{\Xi(e^{-kt})}\|a(t)u(\bxi)v\|^{-\lambda}\,d\bxi + \sum_{i=1}^{\infty} \int_{\Xi(e^{-kt+i})\setminus \Xi(e^{-kt+i-1})}\|a(t)u(\bxi)v\|^{-\lambda}\,d\bxi \\&\ll M(v)^{-\lambda }\sum_{i=0}^{\infty} e^{\lambda (kt-i)}\operatorname{Leb}\big(\Xi(e^{-kt+i})\big)\\&\ll_{\gamma_0} M(v)^{-\lambda} \left(\sum_{0\leq i\leq kt} e^{-(\gamma_0-\lambda)(kt-i)}+\sum_{i>kt}e^{-\lambda(i-kt)}\right)\\&\ll_{\gamma_0} M(v)^{-\lambda}\ll \|v\|^{-\lambda}.
    \end{aligned}
\end{equation*}

\textbf{Subcase 2-2}. $k$ is odd.\\
In this subcase, $\mu_m=1$. Observe that
$\frac{k}{k+1}\lambda<\gamma_k$. Indeed, if $\gamma_k=1$, then
$\lambda<\frac{2R}{n}<1$; otherwise,
$\gamma_k(k+1)=2R>\lambda n\geq\lambda k$. Thus, we may choose $\gamma_0$ so that
$$\frac{k}{k+1}\lambda <\gamma_0<\min\{\gamma_k,\lambda\}.$$
By \eqref{eq:local-two-component-lower-bound}, we have
$$\|a(t)u(\bxi)v\| \geq C^{-1}M(v)e^{-kt}\max\{1, e^{(k+1)t}F(\bxi)\}.$$
In combination with \eqref{eq:sublevelestimateappication}, it follows that
\begin{equation*}
    \begin{aligned}
        \int_{[0,1]^R}\|a(t)u(\bxi)v\|^{-\lambda}\,d\bxi &\leq  \int_{\Xi(e^{-(k+1)t})}\|a(t)u(\bxi)v\|^{-\lambda}\,d\bxi \\&\qquad\qquad+ \sum_{i=1}^{\infty} \int_{\Xi(e^{-(k+1)t+i})\setminus \Xi(e^{-(k+1)t+i-1})}\|a(t)u(\bxi)v\|^{-\lambda}\,d\bxi\\&\ll M(v)^{-\lambda } e^{\lambda k t}\sum_{i=0}^{\infty} e^{-\lambda i}\operatorname{Leb}\big(\Xi(e^{-(k+1)t+i})\big)\\&\ll_{\gamma_0} M(v)^{-\lambda} e^{(\lambda k-\gamma_0(k+1)) t}\sum_{i=0}^{\infty} e^{-(\lambda-\gamma_0)i}\\&\ll_{\gamma_0} M(v)^{-\lambda}\ll \|v\|^{-\lambda}.
    \end{aligned}
\end{equation*}

Combining \textbf{Case 1} and \textbf{Case 2} and using the equivalence of norms completes the proof.

\bigskip

\section{Sharpness examples}\label{section 6}

\begin{example}[Sharpness of $n\geq 4R+1$]
\label{example:sharpness}
Fix $R\geq1$ and a prime $p\equiv3\pmod4$, and put
\begin{equation}\label{eq:qidef}
q_i(\mathbf x)
=x_{4i-3}^2+x_{4i-2}^2
-p(x_{4i-1}^2+x_{4i}^2),
\qquad 1\leq i\leq R.
\end{equation}
This classical block construction is described in \cite{HB13};
see also \cite[p.~1]{C73} for the case $R=3$ and $p=3$.
The system $\mathbf q=(q_1,\ldots,q_R)$ has no non-trivial
common zero over $\mathbb Q_p$, but has the nonsingular
real zero
$
(\sqrt p,0,1,0,\ldots,\sqrt p,0,1,0).$

A standard compactness argument shows that, for sufficiently
large $h$, every system $\mathbf q+p^h\mathbf G$ with
$\mathbb Z_p$-integral coefficients in $\mathbf G$ remains
$\mathbb Q_p$-insoluble.
Indeed, otherwise there would be $h_k\to\infty$ and common
zeros $\mathbf x_k$ of such systems, normalized by
$\|\mathbf x_k\|_p=1$, with
$\|\mathbf q(\mathbf x_k)\|_p\leq p^{-h_k}$.
Then, compactness of
$\mathbb Z_p^{4R}\setminus p\mathbb Z_p^{4R}$
 yields a nonzero common zero of $\mathbf q$.
For this approximation and compactness argument, see
\cite[Section~3, Lemma ~11]{BL65}.

Fix such an $h$ and choose an integral diagonal system
$\mathbf V=(V_1,\ldots,V_R)$ defining a smooth complete
intersection of codimension $R$ in $\mathbb P^{4R-1}$.
Such a choice is possible since the corresponding locus is nonempty and Zariski open in diagonal coefficient space.
For $N\in\mathbb Z$, set
\begin{equation}\label{eq:Qidef}
Q_i^{(N)}=Nq_i+p^hV_i,
\qquad 1\leq i\leq R,
\end{equation}
and let $\mathcal Z_N\subset\mathbb P^{4R-1}$ be their
common zero locus.
Since $\mathcal Z_0$ is a smooth complete intersection
of codimension $R$, Zariski openness gives the same
conclusion for all but finitely many $N$.
For $p\nmid N$, the preceding compactness argument,
applied after dividing the equations by $N$, gives
$\mathcal Z_N(\mathbb Q_p)=\varnothing$, since
$N^{-1}\mathbf V$ has coefficients in $\mathbb Z_p$.
Moreover, coefficientwise over $\mathbb R$,
\[
N^{-1}(Q_1^{(N)},\ldots,Q_R^{(N)})
=\mathbf q+\frac{p^h}{N}\mathbf V
\longrightarrow\mathbf q
\qquad (N\to\infty).
\]
The implicit function theorem at the nonsingular real zero
of $\mathbf q$ displayed above therefore gives
$\mathcal Z_N(\mathbb R)\neq\varnothing$
for all sufficiently large $N$.

Thus, for all sufficiently large positive integers $N$
with $p\nmid N$, $\mathcal Z_N$ is a smooth complete
intersection of codimension $R$ satisfying
$\mathcal Z_N(\mathbb R)\neq\varnothing$ and $
\mathcal Z_N(\mathbb Q_p)=\varnothing.$ 
\end{example}

\begin{example}[Necessity of smoothness]\label{example:smoothness}
Let $n\geq9$ and let $p\equiv3\pmod4$ be a prime. Define
\[
Q_1(\mathbf{x})
=x_1^2+x_2^2-p(x_3^2+x_4^2),
\qquad
Q_2(\mathbf{x})
=\sum_{j=5}^n x_j^2-x_4^2.
\]
Let $\mathcal{Z}\subset \mathbb{P}^{n-1}(\mathbb{Q})$ be the projective variety defined by $Q_1$ and $Q_2$. Since $Q_1$ is irreducible over $\mathbb{C}$ and $Q_2$ is
not a multiple of $Q_1$, the complex variety $\mathcal{Z}_{\mathbb{C}}$ is a complete intersection of codimension $2$.

Observe that the vector $\mathbf{x}_0=(\sqrt p,0,0,1,1,0,\ldots,0)$ is a common real zero. The Jacobian minor in the first and fifth columns equals $4\sqrt p\neq0$, so $[\mathbf{x}_0]$ is a nonsingular real point.

Nevertheless, $\mathcal{Z}_{\mathbb{Q}}=\varnothing$. Indeed, since $-1$ is not a square modulo $p$, reduction of $Q_1(\mathbf{x})=0$ modulo $p$ forces $p\mid x_1,x_2$ for any integer solution. Reduction modulo $p^2$ then forces $p\mid x_3,x_4$. By the infinite descent argument, it follows that $x_1=\cdots=x_4=0$. Then $Q_2(\mathbf{x})=0$ forces all remaining coordinates $x_5,\ldots,x_{n}$ to vanish.
\end{example}

%\bibliographystyle{amsalpha}
%\bibliography{uribib}
\def\cprime{$'$} \def\cprime{$'$} \def\cprime{$'$}
\providecommand{\bysame}{\leavevmode\hbox to3em{\hrulefill}\thinspace}
\providecommand{\MR}{\relax\ifhmode\unskip\space\fi MR }
\providecommand{\MRhref}[2]{%
  \href{http://www.ams.org/mathscinet-getitem?mr=#1}{#2}
}

\end{document}